\documentclass[11pt]{amsart}

\numberwithin{equation}{section}

\newcommand{\lab}{\label}

\newcommand{\ben}{\begin{enumerate}}
\newcommand{\een}{\end{enumerate}}

\newcommand{\bea}{\begin{eqnarray}}
\newcommand{\ba}{\begin{array}}
\newcommand{\bean}{\begin{eqnarray*}}
\newcommand{\ea}{\end{array}}
\newcommand{\eea}{\end{eqnarray}}
\newcommand{\eean}{\end{eqnarray*}}
\newcommand{\beq}{\begin{equation}}
\newcommand{\eeq}{\end{equation}}
\newcommand{\bthm}{\begin{thm}}
\newcommand{\ethm}{\end{thm}}
\newcommand{\blem}{\begin{lem}}
\newcommand{\elem}{\end{lem}}
\newcommand{\bprop}{\begin{prop}}
\newcommand{\eprop}{\end{prop}}
\newcommand{\bcor}{\begin{cor}}
\newcommand{\ecor}{\end{cor}}
\newcommand{\bdfn}{\begin{dfn}}
\newcommand{\edfn}{\end{dfn}}
\newcommand{\brem}{\begin{rem}}
\newcommand{\erem}{\end{rem}}
\newcommand{\bpf}{\begin{proof}}
\newcommand{\epf}{\end{proof}}
\newcommand{\bfact}{\begin{fact}}
\newcommand{\efact}{\end{fact}}

\alph{enumii} \roman{enumiii}

\newtheorem{thm}{Theorem}[section]
\newtheorem{prop}[thm]{Proposition}
\newtheorem{lem}[thm]{Lemma}

\newtheorem{cor}[thm]{Corollary}
\newtheorem{dfn}[thm]{Definition}
\newtheorem{rem}[thm]{Remark}
\newtheorem{fact}[thm]{Fact}

             \def\cB{\mathcal B}       
\def\cH{\mathcal H}             \def\cF{\mathcal F}
           \def\cM{\mathcal M}        
                    \def\cJ{\mathcal J}
\def\cS{\mathcal S}

\def\endpf{\qed}

\def\C{{\mathbb C}}                  \def\oc{\hat \C}

\def\1{1\!\!1}
\def\and{\text{ and }}

\def\dist{\text{{\rm dist}}}

\def\F{{\mathcal F}}

\def\H{\text{{\rm H}}}     \def\HD{\text{{\rm HD}}}

         \def\P{\text{{\rm P}}}

\def\L{{\mathcal L}}

\def\a{\alpha}                \def\b{\beta}             \def\d{\delta}
\def\De{\Delta}               \def\e{\varepsilon}          
\def\g{\gamma}                           \def\l{\lambda}
\def\La{\Lambda}                         
               \def\sg{\sigma}
                          
\def\ka{\kappa}

              \def\bu{\bigcup}
\def\({\bigl(}                \def\){\bigr)}
\def\lt{\left}                \def\rt{\right}

\def\ld{\ldots}               \def\bd{\partial}         \def\^{\tilde}

\def\es{\emptyset}            \def\sms{\setminus}
\def\sbt{\subset}             \def\spt{\supset}
\def\gek{\succeq}

\def\comp{\asymp}
           
\def\sp{\medskip}             \def\fr{\noindent}        
\def\ov{\overline}            

\def\ni{\noindent}

\def\endpf{{\hfill $\square$}}

\def\D{{\mathbb D}}

\newcommand{\amsc}{{\mathbb C}}

\newcommand{\cbar}{\hat{{\mathbb C}} }

\newcommand{\lam}{\lambda}
\newcommand{\ep}{\varepsilon}
\newcommand{\ph}{\varphi}
\newcommand{\al}{\alpha}

\newcommand{\jul}{J(f)}
\newcommand{\fat}{{\mathcal F}_f}
\newcommand{\sing}{sing(f^{-1})}

\newcommand{\exc}{{\mathcal E}_f}

\newcommand{\post}{{\mathcal P}_f}

\newcommand{\pft}{{\mathcal{L}}_t}
\newcommand{\npft}{{\hat{\mathcal{L}}_t}}

\newcommand{\pftl}{{\mathcal{L}_{t,\l }}}
\newcommand{\pftolo}{{\mathcal{L}_{t_0 ,\l _0 }}}

\newcommand{\rad}{{\mathcal J}_r(f)}
\def\fo{f_{\l^0}}
\def\julo{{\mathcal J}(\fo)}
\def\jull{{\mathcal J}(f_\l)}
\newcommand{\radl}{{\mathcal J}_r(f_\l)}
\def\Jgg{{\mathcal J}(g_\g)}

\def\jo{{\mathcal J}_0}

\newcommand{\Jl}{{\mathcal J}_r(f_\l)}
\newcommand{\JPl}{{\mathcal J}_r(f\circ P_\l)}

\begin{document}

\title[]
{ \bf\large {\Large G}eometric Thermodynamical Formalism and Real
Analyticity for Meromorphic Functions of Finite Order}
\date{\today}
\author[\sc Volker MAYER]{\sc Volker MAYER}
\address{Volker Mayer, Universit\'e de Lille I, UFR de Math\'ematiques,
UMR 8524 du CNRS,
59655 Villeneuve d'Ascq Cedex, France}
\email{volker.mayer\@math.univ-lille1.fr\newline \hspace*{0.3cm} Web:
math.univ-lille1.fr/$\sim$mayer}
\author[\sc Mariusz URBA\'NSKI]{\sc Mariusz URBA\'NSKI}
\address{Mariusz Urba\'nski, Department of Mathematics,
 University of North Texas, Denton, TX 76203-1430, USA}
\email{urbanski\@unt.edu\newline \hspace*{0.3cm} Web:
www.math.unt.edu/$\sim$urbanski}
%
%
\thanks{Research of the second author supported in part by the
NSF Grant DMS 0400481.}

\begin{abstract}
Working with well chosen Riemannian metrics and employing
Nevanlinna's theory, we make the thermodynamical formalism work for
a wide class of hyperbolic meromorphic functions of finite order
(including in particular exponential family, elliptic functions,
cosine, tangent and the cosine--root family and also compositions of these functions
with arbitrary polynomials). In
particular, the existence of conformal (Gibbs) measures is established
and then the existence of probability invariant measures
equivalent to conformal measures is proven. As a geometric
consequence of the developed thermodynamic formalism, a version of
Bowen's formula expressing the Hausdorff dimension of the radial
Julia set as the zero of the pressure function and, moreover, the
real analyticity of this dimension, is proved.
\end{abstract}

\maketitle

\section{Introduction}
One of the most fruitful tool in the study of ergodic, stochastic or
geometric properties of a holomorphic dynamical system is the
thermodynamical formalism. We present a completely new uniform
approach that makes this theory available for a very wide class of
meromorphic functions of finite order. The key point is that we
associate to a given meromorphic function $f:\C\to \cbar$ a suitable
Riemannian metric $d\sg= \gamma |dz|$. We then use Nevanlinna's
theory to construct conformal measures for the potentials $-t\log
|f'|_\sg $ and to control the corresponding Perron--Frobenius
operator's. Here
$$
|f'(z)|_\sigma= |f'(z)| \frac{\gamma\circ f(z)}{\gamma(z)}
$$
is the derivative of $f$ with respect to the metric
$d\sigma$. With this tool in hand we obtain then geometric
information about the Julia set $\jul$ and about the radial (or
conical) Julia set
$$\rad =\{z\in \jul \; :\;\; \liminf_{n\to\infty } |f^n(z)| <\infty
\}\; .$$ We now give a precise description of our results.\\

\subsection{Thermodynamical formalism}
 Various versions of thermodynamic formalism and finer fractal
geometry of transcendental entire and meromorphic functions have
been explored since the middle of 90's, and have speeded up since
the year 2000 (see for ex. \cite{bar}, \cite{cs1},
\cite{cs2}\cite{ku1}, \cite{ku2}, \cite{ku3}, \cite{myu1},
\cite{uz0}, \cite{uz1}, \cite{uzdnh}, and especially the survey
article \cite{kuap} touching on most of the results obtained by
now). Some interesting and important classes of functions, including
exponential $\l e^z$ and elliptic, have been fairly well understood.
Essentially all of them were periodic, the methods they were dealt
with broke down in the lack of periodicity, and required to project
the dynamics down onto the appropriate quotient space, either torus
or infinite cylinder. One has actually never completely gone back to
the original phase space, the complex plane $\C$. A nice exception
is the case of critically non-recurrent elliptic functions treated
in \cite{ku2}, where the special but most important potential
$-\HD(J(f))\log|f'|$ was explored in detail. In this paper we
propose an entirely different approach. We do not need periodicity
and we work on the complex plane itself. The main idea, which among
others allows us to abandon periodicity, is that we associate to a
given meromorphic function $f$ a Riemannian conformal metric
$d\sigma = \gamma |dz|$ with respect to which the
Perron-Frobenius-Ruelle (or transfer) operator
\begin{equation}\label{eq pfr}
\pft \ph (w) = \sum _{z\in f^{-1}(w)} |f'(z)|_\sg ^{-t} \ph (z)
\end{equation}
is well defined and has all the required properties that make the
thermodynamical formalism work. Such a good metric can be found for
meromorphic functions $f:\amsc\to \cbar$ that are of finite order
$\rho$ and do satisfy the following growth condition for the
derivative:\\

\

\ni {\it Rapid derivative growth:}
 There are $\al_2 > \max\{0, -\al_1\}$ and $\kappa>0$ such that
 \begin{equation}\label{eq intro growth}
 |f'(z)| \geq  \ka^{-1}(1+|z|^{\al _1})(1+|f(z)|^{\al _2})
 \end{equation}
for all $ z\in \jul \setminus f^{-1}(\infty  )$. Throughout
the entire paper we use the notation
$$
\a=\a_1+\a_2.
$$
This condition is very general and forms our second main idea. It is
comfortable to work with and relatively easy to verify (see
Section~\ref{sec examples}) for a large natural class of functions
which include the entire exponential family $\lam e^z$, certain
other periodic functions ($sin(az+b)$, $\lam\tan(z)$, elliptic
functions...), the cosine-root family $\cos(\sqrt{az+b})$ and the
composition of these functions with arbitrary polynomials. Let us
repeat that in Section~\ref{sec examples} these and more examples are described in
greater detail. The Riemannian metric $\sg$ we are after is
$$
d\sg(z)=(1+|z|^{\a_2})^{-1}|dz|.
$$
Let $(X,m)$ be a probability measure and $T:X\to X$ a measurable
map. Recall that,
given a bounded above non-negative measurable function
$g:X\to[0,+\infty)$, the measure
$m$ is called $g$-conformal provided that
$$
m(T(A))=\int_Agdm
$$
for every measurable subset $A$ of $X$ such that $T|_X$ is
injective. Our third and fourth basic ideas were to revive the old
method of construction of conformal measures from \cite{du1} (which
itself stemmed from the work of Sullivan \cite{sul} and Patterson
\cite{pat}) and to employ results and methods coming from
Nevanlinna's theory. These allowed us to perform the construction of
conformal measures and to get good control of the
Perron-Frobenius-Ruelle operator, resulting in the following key
result of our paper.

\bthm\label{theo main} If $f:\amsc\to \cbar$ is an arbitrary
hyperbolic meromorphic function of finite order $\rho$ that
satisfies the rapid derivative growth condition (\ref{eq intro
growth}), then for every $t> \frac{\rho }{\al}$ the following are
true.
\begin{itemize}
    \item[(1)] The topological pressure $\P(t)=\lim_{n\to\infty}
    \frac{1}{n} \log \pft^n(\1)(w)$ exists and is independent of $w\in J(f)$.
    \item[(2)] There exists a unique $\l|f'|_\sg^t$-conformal
    measure $m_t$ and necessarily $\l=e^{\P(t)}$. Also, there exists
    a unique probability Gibbs state $\mu_t$, i.e.
    $\mu_t$ is $f$-invariant and equivalent to $m_t$. Moreover,
    both measures are ergodic and supported on the radial (or conical) Julia
    set.
    \item[(3)] The density $\psi=d\mu_t/dm_t$ is a continuous and
    bounded function on the Julia set $J(f)$.
\end{itemize}
\ethm

\brem For the existence of $e^{\P(t)}|f'|_\sg^t$-conformal measures
the assumption of hyperbolicity is not needed (see Section \ref{sec
conf m}). \erem

\ni Note that even in the context of exponential functions ($\l e^z$)
and elliptic functions, this result is new since it concerns the map
$f$ itself and not its projection onto infinite cylinder or torus.

An important case in Theorem~\ref{theo main} is when $h$ is a zero
of the pressure
function $t\mapsto \P(t)$. In this situation, the corresponding
measure $m_h$ is $|f'|_\sigma ^h $-conformal (also called simply
$h$-conformal). We will see that such a (unique) zero $h>\rho /\al$
exists provided the function $f$ satisfies a mild growth condition on 
the characteristic function (see (\ref{eq diver})) and, most importantly
the following balanced growth condition.

\

\ni {\it Balanced growth condition:} There are $\al_2 > \max\{0,
-\al_1\}$ and $\kappa>0$ such that
\begin{equation}\label{eq intro growth beta}
 \ka^{-1}(1+|z|^{\al _1})(1+|f(z)|^{\al _2})\leq
 |f'(z)| \leq  \ka (1+|z|^{\al _1})(1+|f(z)|^{\al _2})
\end{equation}
for all finite $ z\in \jul \setminus f^{-1}(\infty  )$.

\

\ni Hyperbolic meromorphic functions of positive and finite order that 
satisfy these conditions are called \it dynamically regular. \rm

\

\subsection{Bowen's formula}

\ni Starting from Section~\ref{gapp} we provide geometric
applications of the key result above and provide, in particular, the
following version of Bowen's formula.

\

\bthm\lab{t1032605} (Bowen's formula) If $f:\C\to\oc$ is a
dynamically regular function,
 then the pressure function $\P(t)$ has a unique zero $h>\rho
/\alpha$ and
$$
\HD(J_r(f))=h \; .
$$
\ethm

\

\ni This  type of formulas has a long and rich history. It has
appeared the first time in the classical Bowen's paper \cite{bow}
and since then has been generalized and adopted to a vast number of
contexts, taking perhaps on the most perfect form in the class of
hyperbolic rational functions. In this class and in many others the
zero of the pressure function is the value of the Hausdorff
dimension of the entire Julia set (which is false for entire
functions \cite{uz0}). By a reasoning, which is by now standard,
Theorem~\ref{t1032605} leads to the following.

\

\bcor\lab{c1032905} With the assumptions of Theorem~\ref{t1032605},
we have $\HD(J_r(f))<2$. \ecor

\

\ni This property applied to the sine or exponential family and
combined with results of McMullen \cite{mcm} (who showed that the
Hausdorff dimension of these functions is always two) gives the
following.

\

\bcor \label{cor 1} If $f$ is any hyperbolic member of the
exponential ($z\mapsto \l e^z$) or the sine ($z\mapsto sin(\a z+\b
)$, $\a\neq 0$) family then the hyperbolic dimension $\HD(J_r(f))$
is strictly less then $\HD(J(f))$. \ecor

Note that such a phenomenon does not exist in the setting of
rational functions. For the exponential family it has been proven in
\cite{uz0}.

\

{\sl Proof of Corollary \ref{c1032905}.} Indeed, by
Theorem~\ref{t1032605} and by Theorem~\ref{theo main} there exists
an $|f'|_\sg^h$-conformal measure for $f$. Suppose to the contrary that
$h=2$. Now the proof is standard (see \cite{uz0} or \cite{my2} for
details): Firstly, using the definition of the set $J_r(f)$, which
gives possibility of taking pull-backs of points lying in a compact
region, and applying Koebe's Distortion Theorem, one shows that the
measure $m_h$ and the $2$-dimensional Lebesgue measure restricted to
$J_r(f)$ are equivalent. Secondly, consider an arbitrary point $z\in
J_r(f)$. As above it has infinitely many pull-backs from a compact
region. Since the Julia set is ``uniformly'' nowhere dense on any
compact part, using Koebe's Distortion Theorem, one easily deduces
that $z$ cannot be a Lebesgue density point of $J_r(f)$. Thus the
Lebesgue measure of $J_r(f)=0$, and this contradiction finishes the
proof. \endpf

\

\subsection{Real analyticity}
Answering the conjecture of D. Sullivan, D. Ruelle in \cite{r} (1982)
gave a proof of the real-analytic dependence of the Hausdorff dimension of the
Julia set for hyperbolic rational maps. More recently, this fact was
extended in \cite{uz1, cs2} to some special families of meromorphic
functions (in particular the exponential family). It was shown that
the variation of the Hausdorff dimension of the radial Julia set
$\rad$ is real-analytic at hyperbolic functions. Note that in the case
of hyperbolic rational functions the Julia and the radial Julia set
coincide. This is no longer true in the meromorphic setting and, as
we have seen in Corollary \ref{cor 1}, there is often a gap between
the \it hyperbolic dimension, \rm i.e. the Hausdorff dimension of
the radial Julia set, and the Hausdorff dimension of the Julia set
itself \cite{uz0}.

\sp

\ni We investigate the variation of the hyperbolic dimension of
meromorphic functions in a very general setting and prove in
particular the following result which contains as special cases the
real analyticity facts established in \cite{uz1} and \cite{cs2}.

\

\bthm\label{1.4} Let $f:\C\to \oc$ be either the sine, tangent,
exponential or the Weierstrass elliptic function and let $f_\l(z)=
f(\l_dz^d +\l_{d-1}z^{d-1} +...+\l_0)$,
$\l=(\l_d,\l_{d-1},...,\l_0)\in \C^*\times\C^{d}$. Then the function
$$
\l \mapsto \HD(\radl)
$$
is real-analytic in a neighbourhood of each parameter $\l^0$ giving
rise to a hyperbolic function $f_{\l^0}$. \ethm

\

\ni This result is an example of an application of the general
Theorem~\ref{1.1} (via Theorem~\ref{t1072005}) that we present now.

The Speiser class
$\cS$ is the set of meromorphic functions $f:\C\to\oc $ that have a
finite set of singular values $\sing$. We will work in the subclass
$\cS_0$ which consists in the functions $f\in \cS$ that have a
strictly positive and finite order $\rho =\rho (f)$ and that are of
divergence type. Fix $\La$, an open subset of $\C^N$, $N\ge 1$. Let
$$
\cM_\La=\{f_\l \in \cS_0\, ;\, \l \in \La\}\; , \;\; \La \subset
\C^N,
$$
be a holomorphic family such that the singular points
$sing(f_\l^{-1} )=\{a_{1,\l }, ..., a_{d,\l } )$ depend continuously
on $\l \in \La$. Consider furthermore $\cH \subset\cS_0$, the set of
hyperbolic functions from $\cS_0$ and put
$$
\cH\cM_\La=\cM_\La\cap\cH.
$$
We say that \ \ni {\it $\cM_\La$ is of bounded deformation} if there
is $M>0$ such that for all $j=1,...,N$ \begin{equation} \label{1.3}
\left|\frac{\partial f_\l(z)}{\partial \l _j}\right| \leq M
|f'_\l(z)| \;\;,\quad \l \in\Lambda \;\; and \;\; z\in \jull .
\end{equation}

\

\ni We also say that $ \cM_\La$ is {\it uniformly balanced}
provided every $f\in  \cM_\La$ satisfies the condition
(\ref{eq intro growth beta}) with some fixed constants $\ka , \al_1 ,
\al_2$.

 \

 \bthm \label{1.1} Suppose $\fo
\in \cH\cM_\Lambda$ is dynamically regular and that $U\subset \La$ is an open neighborhood
of $\l^0$ such that $\cM_U$ is uniformly balanced with $\a_1\geq0$ and of
bounded deformation. Then the map $$\l \mapsto \HD (\radl )$$ is
real-analytic near $\l^0$. \ethm

\section{Generalities}
\label{section2}
 The reader may consult, for example,
\cite{nev'}, \cite{nev} or \cite{hille} for a detailed exposition on
meromorphic functions and \cite{b} for their dynamical aspects. We
collect here the properties of interest for our concerns. The Julia
set of a meromorphic function $f:\amsc\to \cbar$ is denoted by
$\hat{J}(f)$ and the Fatou set by $\fat$. Since we always work in the finite plane
we denote $\jul = \hat{J}(f)\cap \C$. By Picard's theorem, there are
at most two points  $z_0\in \cbar$ that have finite backward orbit
${\mathcal O}^- (z_0)=\bigcup_{n\geq 0}f^{-n}(z_0)$. The set of
these points is the exceptional set $\exc$. In contrast to the
situation of rational maps it may happen that $\exc \subset \jul$.
Iversen's theorem \cite{iv, nev'} asserts that every $z_0\in \exc$
is an asymptotic value. Consequently, $\exc \subset \sing$ the set
of critical and finite asymptotic values. The post-critical set
$\post$ is defined to be the closure in the plane of
$$
\bigcup_{n\geq 0} f^n \big(\sing \setminus f^{-n}(\infty ) \big)\; .
$$
Let us introduce the following definitions.
\
\bdfn\label{defi hyp} A meromorphic function $f$ is called
topologically hyperbolic if
$$
\d(f):={1\over 4}\dist\(\jul,\post\)>0.
$$
and it is called expanding if there is $c>0$ and $\lam >1$ such that
$$
|(f^n)'(z)|\ge c\l ^n \quad for \; all  \;\; z\in \jul \sms
f^{-n}(\infty  )\; .
$$
A topologically hyperbolic and expanding function is called
hyperbolic. \edfn

\ni The Julia set
of a hyperbolic function is never the whole sphere. We thus may and
we do assume that the origin $0\in \fat$ is in the Fatou set
(otherwise it suffices to conjugate the map by a translation). This
means that there exists $T>0$ such that
\begin{equation}\lab{2012705}
D(0,T)\cap J(f)=\es.
\end{equation}
The derivative growth condition (\ref{eq intro growth}) can then be
reformulated in the following more convenient form:\\

{ \noindent \it
 There are $\al_2> 0$, $ \al_1 > -\al_2$ and $\ka>0$ such that
 \begin{equation}\label{eq condition}
 |f'(z)| \geq \ka^{-1} |z|^{\al _1}|f(z)|^{\al _2}\quad for\;
 all  \; \; z\in \jul\setminus f^{-1}(\infty  ) \;.
 \end{equation}
Similarly, the balanced condition (\ref{eq intro growth beta})
becomes
 \begin{equation}\label{eq condition beta}
 \ka^{-1} |z|^{\al _1}|f(z)|^{\al _2}
\leq |f'(z)|
\leq \ka |z|^{\al _1}|f(z)|^{\al _2}\quad for\;
 all  \;\; z\in \jul\setminus f^{-1}(\infty  )
 \end{equation}
and the metric $d\sigma (z) = |z|^{-\a_2}|dz|$.}

\

\ni It is well known that in the context of rational functions
topological hyperbolicity and expanding property are equivalent.
Neither implication is established for transcendental functions.
However, under the rapid derivative growth condition (\ref{eq
condition}) with $\a_1\geq 0$ topological hyperbolicity implies
hyperbolicity.

\

\bprop\lab{p1012805} Every topologically hyperbolic meromorphic
function satisfying the rapid derivative growth condition with
$\a_1\ge 0$ is expanding, and consequently, hyperbolic. \eprop

\

{\sl Proof.} Let us fix $\l \geq 2$
such that $\l \ka^{-1}T^\a \geq 2 $. In view of rapid derivative
growth (\ref{eq condition}) and (\ref{2012705})
\begin{equation}\lab{10050505}
|f'(z)|\ge \ka^{-1}T^\a\;\;\;\; for \; all \;\; z\in J(f)
\end{equation}
 and
\begin{equation}\lab{11050505}
|f'(z)|\ge \l \;\;\;\; for \; all \;\;z\in f^{-1}(J(f)\sms D(0,R))
\end{equation}
provided $R>0$ has been chosen sufficiently large. In addition we
need the following.

\sp\fr {\it Claim:} There exists $ p=p(\l,R)\ge 1$ such that
$$
|(f^n)'(z)|\ge\l\;\;\; for \; all \;\;n\ge p \;\; and \;\; z\in \ov
D(0,R)\cap J(f).
$$

Indeed, suppose on the contrary that there is $R>0$ such
that for some $n_p\to \infty$ and $z_p\in\ov D(0,R)\cap J(f)$ we
have
\begin{equation}\lab{1012805}
|(f^{n_p})'(z_p)|<\l.
\end{equation}
Put $\d=\d(f)$. Then for every $p\ge 1$ there exists a unique holomorphic
branch $f_*^{-n_p}:D\(f^{n_p}(z_p),2\d\)\to\C$ of $f^{-n_p}$ sending
$f^{n_p}(z_p)$ to $z_p$. It follows from ${1\over 4}$-Koebe's Distortion
Theorem and (\ref{1012805}) that
\begin{equation}\lab{2012805}
f_*^{-n_p}\(D\(f^{n_p}(z_p),2\d\)\) \spt D\(z_p,\d /(2\l)\)
\end{equation}
or, equivalently, that $f^{n_p}(D\(z_p,\d /(2\l)\)) \sbt
D\(f^{n_p}(z_p),2\d\)$. Passing to a subsequence we may assume
without loss of generality that the sequence $\{z_p\}_{p=1}^\infty$
converges to a point $z\in \ov D(0,R)\cap J(f)$. Since
$D(\post,2\d)\cap D\(f^{n_p}(z_p),2\d\)=\es$ for every $p\geq 1$, it
follows from Montel's theorem that the family
$\big\{f^{n_p}|_{D(z,(2\l)^{-1}\d)}\big\}_{p=1}^\infty$ is normal,
contrary to the fact that $z\in J(f)$. The claim is proved.

\sp

\ni Let $p=p(\l ,R)\ge 1$ be the number produced by the claim. It
remains to show that
$$
|(f^{2p})'(z)| \geq 2 >1 \quad for \;\; every\;\; z\in J(f).
$$
This formula  holds if $|f^j(z)|>R$ for $j=0,1,...,p$ because
of (\ref{10050505}), (\ref{11050505}) and the choice of $\l$.
If $|z|>R$ but $|f^j(z)|\leq R$ for some $0\le j\leq p$, the
conclusion follows from (\ref{10050505}) and the claim.
\endpf


\

\ni The class of Speiser $\cS$ consists in the functions $f$ that
have a finite set of singular values $\sing$. The classification of
the periodic Fatou components is the same as the one of rational
functions because any map of $\cS$ has no wandering nor Baker
domains \cite{b}. Consequently, if $f\in \cS$ then $f$ is
topologically hyperbolic if and only if the orbit of every singular
value converges to one of the finitely many attracting cycles of
$f$. This last property is stable under perturbation, a fact that we
use in Section \ref{sec 8} and also in the next remark:

\

\bfact\label{2.2} Let $\fo\in \cH$ be a hyperbolic function and
$U\subset \Lambda$ an open neighborhood of $\l^0$ such that, for
every $\l \in U$, $f_\l$ satisfies the balanced growth condition
(\ref{eq condition beta}) with $\ka >0 , \a_1 \geq 0 $ and $\a_2>0$
independent of $\l \in U$. Then, replacing $U$ by some smaller
neighborhood if necessary, all the $f_\l$ satisfy the expanding
property for some $c,\rho$ independent of $\l\in U$. \efact


\

\ni We end this part by giving a more detailed description of the
divergence type functions than the one given in the introduction.
For a meromorphic function $f$ of finite order $\rho$ a theorem of
Borel states that the series
\begin{equation}\label{eq borel sum}
\Sigma(t,w)=\sum_{z\in f^{-1}(w)} |z|^{-t}
\end{equation}
has the exponent of convergence equal to $\rho$ meaning that it
diverges if $t<\rho$ and converges if $t>\rho$. Concerning the
behavior of $\Sigma (t,w)$ in the critical case $t=\rho$ it turns
out that, if $\Sigma (\rho,w)=\infty$ for some $w\in \cbar$, then
this series diverges for all but at most two values $w\in \cbar$
(see Remark~\ref{rem div type}).

\

\begin{dfn}\label{defi div type}
If $\Sigma(\rho,w)=\infty$ for some $w\in \cbar \sms \exc$, then the
function $f$ is said to be of divergence type.
\end{dfn}

\

\ni The symbols $\asymp$ and $\preceq$ will signify through the
whole text that equality respectively inequality holds up to a
multiplicative constant that is independent of the involved
variables.

\section{Functions that satisfy the growth condition}
\label{sec examples}

Here we first explain the meaning of the exponents $\a_1, \a_2$
and then we present various examples that fit into our context.

\subsection{The signification of the exponents $\a_1, \a_2$.}

For entire functions the balanced growth condition (\ref{eq condition beta})
is in fact a condition on the logarithmic derivative of the function.
Indeed, for all known balanced entire functions one has $\a_2=1$ and $\a_1=\rho -1$
with, as usual, $\rho$ being the order of the function. The 
balanced growth condition signifies then that the logarithmic derivative of
the function is of polynomial growth of order $\rho -1$.
The lemma to follow indicates that this is a general fact.
Let $\cB$ be the class of functions that have a bounded singular set
$\sing $. Clearly $\cS\subset \cB$.

\blem \label{lemma VW}
Suppose that $f$ is a balanced entire function of class $\cB$ and of positive finite 
order $\rho$. Then 
$\a_2=1$ and $\a_1=\rho -1$. 
\elem 
 
Notice also that
for all these functions the critical exponent 
$$\Theta = \frac{\rho}{\a}=\frac{\rho}{\a_1+\a_2} =1.$$
We recall that the transfer operator is defined
only for $t>\Theta$ and that this number is the lower
bound given in Bowen's formula. Consequently, entire functions of class $\cB$  
that satisfy the conditions of Bowen's formula do have a hyperbolic dimension
strictly larger then $1$.

\bpf[Proof of Lemma \ref{lemma VW}]
Let $f$ be an entire function as discribed in Lemma \ref{lemma VW}. 
Based on Wiman-Valiron theory, Eremenko \cite{eremenko VW} constructed
an orbit $z_{n+1}=f(z_n)\in I(f)$ with
$$I(f) =\{z\in \C \; , \;\; \lim_{n\to \infty } f^n(z)= \infty\}$$
the escaping set of $f$. Notice that
for entire functions of class $\cB$ we have
$I(f)\subset \jul $ (\cite{el}) and thus $z_n\in\jul $ for every $n\geq 0$.
This orbit has further properties. The important ones for our purpuse
are that, for every $n\geq 0$,
$$|z_{n+1}| \geq \frac{1}{2} m_f(|z_n|)$$
and
$$ \frac{|f'(z_n)|}{|f(z_n)|} = \frac{\nu_f(|z_n|)}{|z_n|}|1+\e _n| \quad where \;\;
|\e_n|< \frac{1}{2} . $$ 
Here $m_f(r)=\max _{|z|=r}|f(z)|$ measures the maximal growth of $f$ and
 the function $\nu_f$ is the so called central index. It follows now easily from the 
growth properties of these functions $m_f$ and $\nu_f$  together with
$$ |f'(z_n)| \asymp |z_n|^{\a_1} |z_{n+1}|^{\a_2}$$
that $\a_2=1$ and $\a_1=\rho -1$.
\epf

For a meromorphic (non entire) function $f$, the exponent $\a_2$ is determined by 
the multiplicity of the function at a pole.
If $q$ is the multiplicity of $f$ at a pole $b$ and if $f$ satisfies the balanced 
growth condition (\ref{eq condition beta})
then 
$$ \a_2 =1+\frac{1}{q} $$
 (cf. the discussion on elliptic functions below).
This looks restrictive in the sense that it forces all poles to have the 
same multiplicity. However this problem can be overcome if one replaces
 the constant $\a_2$ by a function. In order to avoid 
a technical more involved presentation we keep $\a_2$ constant here.

\subsection{Examples}
First of all,
the whole exponential family $f_\lam =\lam \exp (z) $, $\lam \neq
0$, clearly satisfies the growth condition (\ref{eq condition}) with
$\al_1=0$ and $\al_2=1$.  More generally, if $P$ and $Q$ are
arbitrary polynomials, then
$$
f(z) = P(z) \exp (Q(z)) \quad , \;  z\in \amsc \, ,
$$
satisfies (\ref{eq condition}) provided that $|f'|_{|\jul}\geq c
>0$. In this case $\al_1=deg(Q)-1$, $\al_2=1$, the order $\rho=deg
(Q)$ and consequently $\frac{\rho}{\a}=1$.
Assuming still that $|f'|_{|\jul}\geq c >0$ (which holds in
particular for expanding maps), the following functions also satisfy
rapid derivative growth condition (\ref{eq condition}):
\begin{itemize}
    \item[(1)] The sine family: $f(z)= \sin (az+b)$, $a,b\in \amsc$,
    $a\neq 0$.
    \item[(2)] The cosine-root family: $f(z) =\cos(\sqrt{az+b})$ with
    again $a,b\in \amsc$, $a\neq 0$. Note that here
    $\al_1=-\frac{1}{2}$ and $\al_2=1$ which explains that negative
    values of $\al_1$ should be considered in (\ref{eq condition}).
    \item[(3)] Certain solutions of Ricatti differential equations like,
    for example, the tangent family $f(z) = \lam \tan (z) $, $\lam
    \neq 0$, and, more generally, the functions
    $$f(z) =\frac{Ae^{2z^k}+B}{Ce^{2z^k}+D} \quad with \quad AD-BC \neq 0 \; .$$
    The associated differential equations are of the form
    $w'=kz^{k-1}(a+bw+cw^2)$ which explains that here $\al_1=k-1$
    and $\al_2=2$.
    \item[(4)] All elliptic functions.
    \item[(5)] Any composition of one of the above functions with a polynomial.
\end{itemize}
The assertion on elliptic functions deserves some explanation. Let
$f:\amsc \to \cbar$ be a doubly periodic meromorphic function and
let $U=\{z\in\C:|z|>R\} \cup\{\infty\}$, where $R>0$ is chosen
sufficiently large so that:
\begin{itemize}
    \item[a)] every component $V_b$ of $f^{-1}(U)$ is a bounded
    topological disc, and
    \item[b)] there is $\kappa >0$ such that for every pole $b$ and
    any $z\in V_b \setminus \{b\}$ we have $|f'(z)| \asymp
    |f(z)|^{1+\frac{1}{q_b}}$ where $q_b$ is the multiplicity of the
    pole $b$.
\end{itemize}
From the periodicity of $f$ and the assumption $|f'|_{|\jul}\geq c
>0$ easily follows now that $f$ satisfies (\ref{eq condition}) with
$\al_1=0$ and
$$
\al_2 =\inf \Big\{1+\frac{1}{q_b}:b\in f^{-1}(\infty)\Big\}\quad .
$$
More generally, the preceding discussion shows that for any function
$f$ that has at least one pole one always has $$ \al_2 \leq \inf
\Big\{1+\frac{1}{q_b}:b\in f^{-1}(\infty)\Big\} .$$

\

\ni The stronger balanced growth condition (\ref{eq condition beta})
is also satisfied by an elliptic function provided all its poles
have the same order. General elliptic and meromorphic functions with
poles of different order can not be of balanced
growth. As already mentioned, this problem can be overcome if one allows
$a_2$ being variable; we clarify this situation in a forthcoming paper
(were also other classes of balanced functions will be given).

\

\ni Uniform balanced growth is verified by various families. Here are
some examples.

\

\blem \label{lemm unif balanced}
Let $f:\C\to \oc$ be either the
sine, tangent, exponential or the Weierstrass elliptic function and
let $f_\l(z)= f(\l_dz^d +\l_{d-1}z^{d-1} +...+\l_0)$,
$\l=(\l_d,\l_{d-1},...,\l_0)\in \C^*\times\C^{d}$. Suppose $\l^0$ is
a parameter such that $f_{\l^0}$ is topologically hyperbolic. Then
there is a neighbourhood $U$ of $\l^0$ such that $\cM_U=\{f_\l \; ;
\l\in U\}$ is uniformly balanced.
\elem
\bpf All the functions $f$ mentioned have only finitely many
singular values, they are in the Speiser class. The function
$f_{\l^0}$ being in addition topologically hyperbolic, its singular
values are attracted by attracting cycles. As we already remarked in
the previous section, this is a stable property in the sense that
there is a neighbourhood $U$ of $\l^0$ such that all the functions
of $\cM_U=\{f_\l \; ; \l\in U\}$ have the same property. In
particular, no critical point of $f_\l$ is in $J(f_\l )$.
The function $f$ satisfies a differential equation of the form
$$ (f')^p = Q\circ f $$
with $Q$ a polynomial whose zeros are contained in $sing(f^{-1})$.
For example, in the case when $f$ is the Weierstrass elliptic
function then
$$ (f')^2 = 4 (f-e_1)(f-e_2)(f-e_3)$$
with $e_1,e_2,e_3$ the critical values of $f$.
Let $\l \in U$ and denote $P_\l (z) = \l_dz^d +\l_{d-1}z^{d-1}
+...+\l_0$. Since
$$(f_\l ')^p = (f'\circ P_\l \,P_\l ')^p = Q\circ f_\l (P_\l
')^p\; $$ and $f_\l (z)\neq 0$ for all $z\in J(f_\l )$, the
polynomials $P_\l '$ and $Q$ do not have any zero in $J(f_\l)$.
Consequently
$$ |P_\l ' (z)|\asymp |z|^{d-1}     \quad and \quad |Q(z)| \asymp
|z|^q \quad on \quad J(f_\l )$$ with $q=deg(Q)$. Moreover,
restricting $U$ if necessary, the involved constants can be chosen
to be independent of $\l \in U$. Therefore,
$$|f_\l '(z) | \asymp |f_\l(z)|^{q\over p} |z|^{d-1}$$
for $z\in J(f_\l)$ and $\l \in U$. We verified the uniform balanced
growth condition with $\a_1= d-1$ and $\a_2=\frac{q}{p}$ depending
on the choice of $f$. In the case of the Weierstrass elliptic
function one has $\a_2= 3/2$. \epf

\section{Growth condition and cohomological transfer operator}
\ni For exponential or elliptic functions one can use the periodicity to
project the map onto the quotient space (torus or cylinder). This
idea recently lead to many new results (see \cite{kuap} and the
reference therein). Here we replace
the quotient spaces by metric spaces $(\amsc,d\sigma )$ which are
much more flexible. The first and essential problem however is to
find the right natural metric for a given meromorphic function. We
will describe now how this can be done for meromorphic functions of
finite order that satisfy the rapid derivative growth condition.
Recall that we work with the metric
$$
d\sigma(z) =|z|^{-\al _2}|dz|
$$
and we set $\al =\al_1+\al_2$. The derivative of a
function $f:\amsc \to \cbar$ with respect to this metric is given at
a point $z\in\C$ by the formula
$$
|f'(z)|_\sigma
= \frac{d\sigma(f (z))}{d\sigma (z)}
=|f'(z)|\frac{|f(z)|^{-\a_2}}{|z|^{-\a_2}}
=|f'(z)||z|^{\a_2}|f(z)|^{-\a_2}.
$$
We will now see that this is the right choice of the metric in order
for the associated transfer operator $\L_t$ (with the potential
$-t\log|f'|_\sigma $) to act continuously on the Banach space
$C_b(J(f))$ of bounded continuous functions on $J(f)$. Indeed
$$
\aligned \pft \ph (w) &= \sum_{z\in f^{-1}(w)} |f'(z)|_\sigma
^{-t}\ph (z)
= \sum_{z\in f^{-1}(w)} |f'(z)|^{-t}|z|^{-\a_2t}|f(z)|^{\a_2t}\ph (z) \\
&=|w|^{\a_2t}\sum_{z\in f^{-1}(w)} |f'(z)|^{-t}|z|^{-\a_2t}\ph (z).
\endaligned
$$
So, if $f$ satisfies (\ref{eq condition}), then
\begin{equation}\label{eq 3}
\pft \1 (w)
=\sum _{z\in f^{-1}(w)}|f'(z)|_\sigma^{-t}
\leq \ka^t \sum _{z\in f^{-1}(w)}|z|^{-\a t}.
\end{equation}
Now, assume that $f$ is of finite order $\rho$. Then, as we noted in
the introduction, a theorem of Borel states
that the last series has the exponent of convergence equal to $\rho$
for all but at most two points (the points from $\exc$). Assume that
$\exc$ is disjoint from the Julia set $J(f)$; this is for example true
if $f$ is topologically hyperbolic. What we need is the
uniform convergence of the last series in (\ref{eq 3}) in order to
secure continuity of the
operator $\pft$ on the Banach space $C_b(J(f))$ of bounded continuous functions
endowed with the standard supremum norm. More
precisely, we need to know that, for a given $t>\rho/\a$, there
is $M_t>0$ such that
\begin{equation}\label{eq 6}
\pft \1(w) \leq M_t \quad for \; all \;\; w\in \jul \; .
\end{equation}
It turns out that under our assumptions this is always true:

\

\bthm\label{prop 4} Assume that $f:\C\to\oc$ is a finite order
hyperbolic meromorphic function of rapid derivative growth. Then for
every $t> \rho/\a$, the transfer operator $\L_t$ is well defined and
acts continuously on the Banach space $C_b(J(f))$. \ethm

\

\ni The rest of this section is devoted to the proof of this "Uniform
Borel Theorem". Our proof relies on Nevanlinna's theory of value distribution.
The reader can find in the modern monograph \cite{cy} a complete account
on this topic. 
Let $f$ be meromorphic of finite order $\rho$ and
let $u> \rho$. We are interested in the dependence of the following
series on $a$:
$$
\Sigma(u,a)=\sum _{f(z) = a} \frac{1}{|z|^u} \; .
$$
Borel's theorem states that this series converges for every
non-exceptional value $a\in \cbar$. But is this convergence uniform?
To see this we investigate the error terms in the proof of Borel's
theorem as given in \cite[p. 265]{nev'} or \cite[p. 261]{nev}. In
order to do this we use again the fact that $0\in \fat$. In the
following we use the standard notations of Nevanlinna theory. For
example, $n(t ,a)$ is the number of $a$-points of modulus at most
$t$, $N(r,a)$ is defined by $dN(r,a)= n(r,a)/r$ and $T(r)$ is the
characteristic of $f$ (more precisely the Ahlfors-Shimizu version of
it; these two different definitions of the characteristic function only
differ by a bounded amount). The first main theorem (FMT) of
Nevanlinna yields the following for our situation:

\

\bcor[\bf of FMT] \label{fmt}
There is $\Xi>0$ such that $N(r,a)\leq T(r) +\Xi$ for all $a\in \jul $.
\ecor

{\sl Proof.} FMT as stated in \cite{eremenko} or in \cite[p.
216]{hille} yields
$$N(r,a) \leq T(r) +m(0,a) \quad for \; all \; r>0 \; and \; a\in
\cbar $$ with $m(0,a) =-\log [f(0), a]$ and where $[a,b]$ denotes
the chordal distance on the Riemann sphere (with in particular
$[a,b]\leq 1$ for all $a,b\in \cbar$). Since $f(0) \in \fat$, there
is $\tau >0$ such that $[a, f(0)] >\tau$ for all $a\in \jul$. It
follows that the error term is bounded by
$$
0\leq m(0,a) \leq -\log \tau =\Xi \quad for \; all \; \;a\in \jul \;
.
$$
\endpf

\

\ni From the second main theorem (SMT) of Nevanlinna we need the
following version which is from \cite[p. 257]{nev'} (\cite[p.
255]{nev} or again \cite{hille}) and which is valid only since $f$
is supposed to be of finite order.

\

\bcor[\bf of SMT] \label{smt} Let $a_1,a_2,a_3\in \cbar$ be distinct
points. Then
$$
T(r)
\leq \sum_{j=1}^3 N(r,a_j) + S(r) \quad for \; every \;\; r>0\;\;
  with \;\;S(r)
= {\mathcal O} (\log (r)).
$$
\ecor

\

\ni The error term $S(r)$ has been studied in detail and sharp estimates
are known. The following results from Hinkkanen's paper \cite{hinkkanen}
and also from Cherry-Ye's book \cite{cy}.

\blem \label{lemma sharp smt}
Let $f$ be a hyperbolic meromorphic function of finite order $\rho$ that is 
normalized such that $0\in D(0,T)\subset \fat $, $f(0)\notin\{0,\infty \}$
and $f'(0)\neq 0$. Then, for every $\Delta <T/4$, there exists 
$C_1=C_1(\Delta )>0$ and $C_2>0$ such that
$$ 4N(R+\Delta ,a) \geq T(R) -(3\rho +1)\log R -C_1 - C_2 \log |a|  $$
for every $a\in \jul$ and every $R>T$.
\elem

\

\bpf
Since $f$ is expanding there is $c>0$ such that $|f'(z)|\geq c >0$ for all $z\in\jul$.
Let $0<\Delta ' < \min \{\delta (f) , T \}$ such that
$\Delta =2K\Delta'/c <T/4$ where $K$ is an appropriate Koebe distortion constant.
Consider then $a\in \jul$ and $a'\in D(a ,\Delta')$.  Since all the inverse branches 
of $f$ are well defined on $D(a ,2\Delta' )$ we have
$$ n(r+\Delta ,a) \geq n(r,a') \quad , \quad r>0.$$
Consequently
\begin{eqnarray*}
 N(R,a' ) & =& \int_0^R \frac{n(r,a')}{r}\, dr \leq \int_0^R \frac{n(r+\Delta ,a)}{r}\, dr\\
&=& \int_{\Delta}^{R+\Delta} \frac{n(t ,a)}{t}\frac{t}{t-\Delta}\, dt
\leq \frac{T}{T-\Delta}\int_T^{R+\Delta} \frac{n(t ,a)}{t}\, dt\\
&\leq& \frac{4}{3} N(R+\Delta ,a) \;\; for\;  every \;\,R>T.
\end{eqnarray*}
Choose now $a_1,a_2,a_3\in D(a,\Delta')$, any points that satisfy 
$|a_i-a_j|\geq \Delta'/3$ for all $i\neq j$. It follows then from the 
sharp form of SMT given in \cite{hinkkanen}, the fact that $f$ is of finite order, 
along with the normalisations stated in the lemma that
\begin{eqnarray*}
4N(R+\Delta , a) &\geq& \sum_{i=1}^3 N(R,a_i) \geq T(R)-S(R, a_1, a_2, a_3 ) \\
&\geq & T(R)-(3\rho +1)\log R -C_1(\Delta ) - C_2 \log |a| 
\end{eqnarray*}
for every $a\in\jul $ and for all $R>T$.
\epf

\

\ni We can now show the following uniform version of Borel's theorem
which implies Theorem~\ref{prop 4}.

\

\bprop \label{prop borel effectif} Let $f$ be meromorphic of finite
order $\rho$ and suppose that $0\in \fat$. Then, for every $u
>\rho$, there is $M_u>0$ such that
$$ \Sigma(u,a)
=\sum _{f(z)= a} \frac{1}{|z|^u}
\leq M_u\quad for \; all \;\; a\in \jul
\;.
$$
\eprop

{\sl Proof.} Recall that $\jul \cap \ov D(0,T) =\emptyset$.
Then $n(T,a) = N(T, a) = 0$, for all $a\in \jul$, and
by the definition of the Riemann-Stieltjes integral, integration by
parts and the fact that $\lim_{r\to\infty}\frac{n(r,a)}{r^{u}}
=0$, we get that
$$
\Sigma(u,a)= \int_T^\infty\frac{d\, n(t ,a)}{t^u} = u
\int_T^\infty \frac{n(t,a)}{t^{u+1}} \, dt\; .
$$ In the
same way
$$ \int_T^\infty \frac{n(t,a)}{t^{u+1}} \, dt=u
\int_T^\infty \frac{N(t,a)}{t^{u+1}} \, dt \; .
$$
Putting both equations together, we get
\begin{equation}\label{eq 0}
\Sigma(u,a)= u^2 \int_T^\infty \frac{N(t,a)}{t^{u+1}} \,
dt\; .
\end{equation}
Now we proceed like in the proof of Borel's theorem as stated in
\cite[p. 265]{nev'} or in \cite[p. 261]{nev}: let $a_1,a_2,a_3$ be
three different points of $\jul$ and let $a\in \jul $ be any point.
Then it follows from FMT and SMT as stated above that, for every
$t>T$,
\begin{equation}\label{eq 1}
N(t,a) -\Xi \leq T(t) \leq N(t,a_1)+N(t,a_2)+N(t,a_3) + S(t)\; .
\end{equation}
Dividing this relation by $t^{u+1}$ and integrating with respect
to $t$ gives
$$
\int_T^\infty \frac{N(t,a)}{t^{u+1}} \, dt \leq \sum
_{j=1}^3\int_{T}^\infty \frac{ N(t,a_j)}{t^{u+1}} \, dt + A_u\; .
$$
Here we used the fact that $S(r) ={\mathcal O}(\log (r))$, which
implies that $\int_T^\infty \frac{S(t)}{t^{u+1}} \, dt = A_u<
\infty$. Together with (\ref{eq 0}) we finally have
\begin{equation}\lab{1051205}
\Sigma(u,a) \leq \Sigma(u,a_1)+\Sigma(u,a_2)+\Sigma(u ,a_3) +u^2 A_u
\end{equation}
for every $a\in J(f)$.
\endpf

\

\brem\label{rem div type} If the order $\rho>0$, then the above
proof shows that $\Sigma (\rho , b)=\infty$ for all but at most two
values $b\in \cbar$ provided $\Sigma (\rho , a) =\infty$ for some
$a\in \cbar$. This property trivially also holds if $\rho=0$.
Note that Koebe's distortion theorem and hyperbolicity yield that
the two exeptional values for this property cannot be in $\jul$. Therefore
$\Sigma (\rho , a)=\infty$ for all or none $a\in \jul$ and this occurs if and 
only if
$$ \int \frac{T(r)}{r^\rho}\, dr =\infty \, .$$
\erem

\

\section{Construction of conformal measures}\label{sec conf m}

 \ni Further properties of transfer operators $\L_t$ rely on
the existence of conformal measures. Define now the topological pressure as
follows.
\begin{equation}\label{eq 5}
\P(t)=\P(t,x) =\limsup_{n\to\infty } \frac{1}{n} \log \pft^n \1 (x).
\end{equation}
Note that because of hyperbolicity of the function $f$, Koebe's
Distortion Theorem and density in $J(f)$ of the full backward orbit
of any point in $J(f)$, the number $\P(t)=\P(t,x)$ is independent
of $x\in J(f)$.
We recall that $m_t$ is called $e^{\P(t)}|f'|_\sigma^t
$-conformal if $\frac{dm_t \circ f }{ dm_t} =e^{\P(t)}|f'|_\sigma^t$
or, equivalently, if $m_t$ is an eigenmeasure of the adjoint $\pft^*$
of the transfer operator $\pft$ with eigenvalue $e^{\P(t)}$. Note that
then the measure $m_t^e$, the Euclidean version of $m_t$, defined by
the requirement that $dm_t^e(z)=
|z|^{\a_2t}dm_t(z)$ is $e^{\P(t)}|f'|^t$-conformal. If $\P(t)=0$, then
these measures are called $t$-conformal.
In \cite{sul} Sullivan has proved that every rational function
admits a probability conformal measure. As it is shown below, in the
case of meromorphic functions the situation is not that far apart.
All what you need for the existence of an
$e^{\P(t)}|f'|_\sigma^t$-conformal measure is the rapid
derivative growth; no hyperbolicity is necessary.
\footnote{If $f$ is not hyperbolic then $\fat =\emptyset$ may occur
and our method does not work. But then the Lebesgue
measure is $2$-conformal.}
 We present here a
very general construction of conformal measures.

\

\bthm\label{theo conf meas}
If $f:\amsc\to \cbar$ is a meromorphic
function of finite order with non-empty Fatou set satisfying the
rapid derivative growth condition, then for every
$t>\rho/\a$ there exists a Borel probability
$e^{\P(t)}|f'|_\sg^t$-conformal measure $m_t$ on $\jul$.
\ethm

\

\ni The rest of this section is devoted to the proof of Theorem
\ref{theo conf meas}. First of all, changing the system of coordinates
by translation, we may assume without loss of generality that $0\notin
J(f)$. Fix $x\in J(f)$. Observe that the transition parameter for
the series
$$
\Sigma_s =\sum_{n=1}^\infty e^{-ns} \pft ^n \1 (x)
$$
is the topological pressure $\P(t)$. In other words, $\Sigma _s
=+\infty$ for $s<\P(t)$ and
$\Sigma_s <\infty$ for $s>\P(t)$. We assume that we are in the
divergence case, e.g. $\Sigma _{\P(t)} = \infty$. For the convergence
type situation the usual modifications have to be done (see
\cite{du1} for details). For $s>\P(t)$, put
$$
\nu_s =\frac{1}{\Sigma _s} \sum_{n=1}^\infty e^{-ns} (\pft ^n
)^* \delta _x \; .
$$
The following lemma follows immediately from definitions.

\blem\lab{prop 1}
The following properties hold:
\begin{enumerate}
    \item For every $\ph\in {\mathcal C}_b (\amsc )$ we have
$$
\int\ph d\nu_s
=\frac{1}{\Sigma _s} \sum_{n=1}^\infty e^{-ns}\int\pft ^n \ph d\d_x
=\frac{1}{\Sigma _s} \sum_{n=1}^\infty e^{-ns}\pft ^n \ph (x) \; .
$$
    \item  \hspace{2.1cm} $\nu_s$ is a probability measure.
\end{enumerate}
$$ (3) \qquad \qquad  \qquad \frac{1}{ e^s} \pft ^* \nu_s
   = \frac{1}{\Sigma _s} \sum_{n=1}^\infty e^{-(n+1)s}
    (\pft^{n+1})^* \delta _x = \nu_s - \frac{1}{\Sigma _s}
    \frac{\pft^* \delta _x}{ e^s} \; . \qquad \qquad \qquad
$$
\elem

\

\ni The key ingredient of the proof of Theorem~\ref{theo conf meas}
is to show that the family $(\nu_s)_{s>P(t)}$ of Borel probability
 measures on $\C$ is tight and then
to apply Prokhorov's Theorem. In order to accomplish this we put
$$
U_R=\{z\in\C:|z|>R\}
$$
and start with the following observation.

\

\blem\label{lemm 1}
For every $t>\rho/\a$ there is $C=C(t)>0$ such that
$$
\pft (\1_{U_R}) (y) \leq \frac{C}{R^{\al\gamma}} \;\; for \;every
\;\;y\in \jul ,
$$
where $\gamma=\frac{t-\rho/\a}{2}$.
\elem
{\sl Proof.} From the rapid derivative growth condition
(\ref{eq condition}) and
Proposition~\ref{prop borel effectif}, similarly as (\ref{eq 3}), we
get for every $y\in J(f)$ that
\begin{eqnarray*}
\pft (\1_{U_R}) (y) &
=& \sum_{z\in f^{-1}(y)\cap U_R}|f'(z)|_\sg^{-t}
\leq \ka^t \sum_{z\in f^{-1} (y) \cap U_R} |z|^{-\alpha t}\\
&\leq &\frac{\ka^t}{R^{\alpha \gamma}} \sum_{z\in f^{-1} (y)} |z|
^{-(\rho +\a \gamma)}
\le \frac{\ka^tM_{\rho +\a \gamma}}{R^{\alpha \gamma}}.
\end{eqnarray*}
\hfill\endpf

\

\ni Now we are ready to prove the tightness we have already
announced. We recall that this means that
\begin{equation*}
\forall \ep >0 \;\; \exists R>0 \;\; such \; that \;\; \nu_s (U_R)
\leq \ep  \;\; for \; all \; s>\P(t)\; .
\end{equation*}

\

\blem\label{prop 2} The family $(\nu_s)_{s>\P(t)}$ of Borel
probability measures on $\C$ is tight and, more precisely, there is
$L>0$ and $\delta >0$ such that
\begin{equation*}
 \nu_s (U_R)
\leq L R^{-\delta}   \;\; for \; all \;R>0\; and\;  s>\P(t)\; .
\end{equation*}
\elem {\sl Proof.} The first observation is that
\begin{eqnarray*}
\pft^{n+1}(\1_{U_R}) (x)&
= &\sum_{y\in f^{-n}(x)} \sum_{z\in f^{-1}(y)\cap U_R}
   \Big(|f'(z)|_\sigma \,|(f^n)'(y)|_\sigma\Big)^{-t}\\
&=&\sum_{y\in f^{-n}(x)}|(f^n)'(y)|_\sigma^{-t}\pft(\1_{U_R}) (y)
\leq \frac{C}{R^{\al \gamma}} \pft ^n \1 (x).
\end{eqnarray*}
where the last inequality follows from Lemma~\ref{lemm 1}.
Therefore, for every $s> \P(t)$, we get that
\begin{eqnarray*}
\nu_s (U_R ) &=& \frac{1}{\Sigma _s} \sum_{n=1}^\infty e^{-ns}
 \pft^n (\1_{U_R}) (x) \leq \frac{C}{R^{\al \gamma}} \frac{1}{\Sigma _s}
\sum_{n=1}^\infty e^{-ns} \pft ^{n-1}
\1 (x)\\
&=& \frac{C}{R^{\al \gamma}} \frac{1}{ e^s}\frac{1}{ \Sigma _s }
\Big( 1+ \sum_{n=1}^\infty e^{-ns} \pft ^{n} \1 (x)\Big)\leq
\frac{2C}{e^{P(t)}} \frac{1}{R^{\al \gamma}}.
\end{eqnarray*}
This shows Lemma \ref{prop 2} and the tightness of the family
$(\nu_s)_{s>\P(t)}$. \hfill $\square$

\

\ni Now, choose a sequence $\{s_j\}_{j=1}^\infty$, $s_j>\P(t)$, converging
down to $\P(t)$. In view of Prokhorov's Theorem and Lemma~\ref{prop 2},
passing to a subsequence, we may assume without loss of generality that
the sequence $\{\nu_{s_j}\}_{j=1}^\infty$ converges weakly to a Borel probability
measure $m_t$ on $J(f)$. It follows from Lemma~\ref{prop 1} and the divergence
property of $\Sigma _s$ that $\pft^* m_t = e^{P(t)} m_t$. The proof of
Theorem~\ref{theo conf meas} is complete.

\

\section{Gibbs states}
\ni We now complete the proof of Theorem~\ref{theo main}. The first
observation is that one can have a better estimate than
(\ref{eq 3}) in diminishing $\al_2$ slightly. Suppose that the
derivative of $f$ satisfies the growth condition (\ref{eq
condition}) with $\al _2 ' =\al_2+\ep$, $\ep >0$, instead of
$\al_2$. Then
$$
|f'(z)|_\sg = \frac{|z|^{\al_2}}{|f(z)|^{\al_2}} |f'(z)| \geq
\frac{1}{\ka} |z|^{\al} |f(z)|^{\ep} \; , \quad z\in \jul\, ,
$$
which, along with Proposition~\ref{prop borel effectif},
leads to the following important estimate of the transfer operator.
For each $t>\rho/a$,
\begin{equation}\label{eq 4}
\pft \1 (w)
\leq \frac{\ka^t}{|w|^{t\ep}} \sum_{z\in f^{-1}(w)} |z|^{-t\al}
\leq \frac{\ka^tM_{\a t}}{|w|^{t\ep} } \quad for \; all \;\; w \in
\jul \, .
\end{equation}
An immediate advantage of this estimate is the following.

\

\blem\label{lemm 2}
We have $\lim_{w\to\infty} \pft \1 (w) =0$.
\elem

\

\ni The last ingredient we need in this section is the following
straightforward consequence of Proposition~\ref{p1012805}, Koebe's
Distortion Theorem, and the fact that $0\notin J(f)$.

\

\blem\lab{sgbdt} For every hyperbolic meromorphic function
$f:\C\to\oc$ satisfying the rapid derivative growth condition there
exists a constant $K_\sg\ge 1$, called $\sg$-adjusted Koebe
constant, such that if $R>0$ is sufficiently small, then for every
integer $n\ge 0$, every $w\in J(f)$, every $z\in f^{-n}(w)$ and all
$x,y\in D_\sg(w,R|w|^{-\a_2})\cup D(w,R)$ , we have that
\begin{equation}\lab{1032805}
K_\sg^{-1}\le {|(f_z^{-n})'(y)|_\sg\over |(f_z^{-n})'(x)|_\sg} \le K_\sg.
\end{equation}
\elem

\

\ni As an immediate consequence of this lemma and Montel's theorem,
which implies that for every open set $U$ intersecting the Julia set $J(f)$
and every point $z\in J(f)$ there exists $n\ge 0$ such that
$U\cap f^{-n}(z)\ne\es$, we conclude that the topological pressure
$$
\P(t)=\lim_{n\to\infty}\frac{1}{n}\log \pft^n(\1)(w)
$$
exists and is independent of $w\in J(f)$. From these two
lemmas above and the existence of conformal measures
(Theorem~\ref{theo conf meas}) one gets, following the arguments
from formula (3.6) through Lemma~3.6 of \cite{uz1}, the following
uniform estimates for the normalized transfer operator
$$
\npft = e^{-\P(t)} \pft.
$$

\

\bprop\label{prop 3}
There exists $L>0$ and, for every $R>0$, there exists $l_R>0$ such that
$$
l_R \leq  \npft ^n \1(w) \leq L
$$
for all $ n\geq 1$ and all $w\in \jul\cap D(0,R)$. \eprop

\

\ni This allows us to construct an everywhere positive, decreasing to
zero at infinity, fixed point $\psi$ of the normalized transfer operator
$\npft$ by putting
$$
\psi= \tilde{\psi_t}\Big/ \int \tilde {\psi_t} dm_t \quad with \quad
\tilde{\psi_t}(z) = \liminf _{n\to\infty} \frac{1}{n} \sum _{k=1}^n \npft
^k \1 (z) \; , \; z \in \jul
$$
The Borel probability measure $\mu_t=\psi_tm_t$ is obviously $f$-invariant and
equivalent to $m_t$. Repeating the appropriate reasonings from \cite{uz1}
or \cite{myu1}, the proof of Theorem~\ref{theo main} follows.

\

\section{Geometric applications}\lab{gapp}

\ni In the rest of the paper we derive several geometric
consequences from the dynamical results proven in the previous
sections. Our primary goal is to complete the proof of
Theorem~\ref{t1032605} (Bowen's formula). For this part we
strengthen our assumptions and assume
throughout the whole rest of the paper that $f$ is dynamically regular.

\

\begin{dfn} \label{dyn reg} The meromorphic function
 $f:\C\to\oc$ is called dynamically regular if it is hyperbolic,
 of positive and finite order $\rho$,
satisfies the balanced derivative growth (condition (\ref{eq condition
beta})) and if it is of divergence type. In the case $f$ is entire we assume
instead of the divergence type assumption that, for any $A,B>0$, there exists $R>1$
such that
\begin{equation}\label{eq diver}
\int_{\log R}^R \frac{T(r)}{r^{\rho +1}}\, dr - B \left(\log R\right)^{1-\rho} \geq A.
\end{equation}
\end{dfn}

\

\ni The divergence type assumption and also (\ref{eq diver}) do hold in particular
if 
$$ \liminf_{r\to\infty} \frac{T(r)}{r^\rho}>0.$$
This last condition is satisfied by all the examples given in Section \ref{sec examples}.

\

\ni In order to bring up geometric consequences, we
need some information about the shape of the graph of the pressure
function.

\

\bprop\lab{p1051205} If $f:\C\to\oc$ is dynamically regular, then the
following hold.
\begin{itemize}
\item[(a)] The function $t\mapsto \P(t)$, $t>\rho/\a$, is convex and,
consequently, continuous.
\item[(b)] The function $t\mapsto \P(t)$, $t>\rho/\a$, is strictly decreasing.
\item[(c)] $\lim_{t\to+\infty}\P(t)=-\infty$.
\item[(d)] $\lim_{t\to(\rho/\a)^+}\P(t)>0 $.
\end{itemize}
\eprop

\

{\sl Proof.} Convexity of the pressure function $\P(t)$ follows
immediately from its definition and H\"older's inequality. So, item
(a) is proved. Items (b) and (c) are straightforward consequences of
the expanding property.

 It remains to show item (d). If $f$ has a pole $b$ of multiplicity $q$ 
then, since $f$ is balanced, $\a_2=1+1/q$ (see the discussion on elliptic functions
in Section \ref{sec examples}). The result in \cite{my} shows then that
$$ \P(\rho / \a ) \geq 0$$
with strict inequality if $f$ is of divergence type (see Remark 3.2 of \cite{my}).
 
So, let finally $f$ be entire. 
 Notice first that, with the balanced growth condition, the calculations
leading to (\ref{eq 3}) give the following lower estimate.
\begin{equation}\lab{4051305}
\pft \1 (w)
=\sum _{z\in f^{-1}(w)}|f'(z)|_\sigma^{-t}
\ge \ka^{-t} \sum _{z\in f^{-1}(w)}|z|^{-\a t} \quad , \; w\in J(f) ,
\end{equation}
for all $t>\rho/\a$. Denote now for every $R>0$
$$ \Sigma ^R (u,a) = \sum_{z\in f^{-1}(a)\cap D(0,R)} |z|^{-u} \quad , \,\, a\in \jul \; and \; u\geq \rho.$$
In order to proof $\P(\rho/ \a )>0$ it suffices to show that for a given $A>0$ there exists 
$R>T$ such that
$$ \Sigma ^R (\rho ,a) \geq  A  \quad for \; all \;\; a\in \jul \cap D(0,R).$$
Let $R>T$ and let $a\in D(0,R)\cap \jul$. We get precisely in the same way as in (\ref{eq 0}) that
$$\Sigma ^R (\rho ,a) \geq \rho^2 \int_0^R \frac{N(t,a)}{t^{\rho +1}}\, dt 
\geq \int _{\log |a|}^R \frac{N(t,a)}{t^{\rho +1}}\, dt.$$
From the sharp form of the SMT (Lemma \ref{lemma sharp smt}), it follows that
\begin{eqnarray*}
 \Sigma ^R (\rho ,a)& \geq &\rho^2 \int_{\log |a|-\Delta }^{R-\Delta } \frac{N(r+\Delta , a)}{r^{\rho +1}}
\left(\frac{r}{r+\Delta } \right)^{\rho +1} dr \asymp
 \int_{\log |a|-\Delta }^{R-\Delta } \frac{4N(r+\Delta , a)}{r^{\rho +1}} dr\\
&\geq & \int_{\log |a|-\Delta }^{R-\Delta } \frac{T(r)}{r^{\rho +1}} dr -C_1 -C_2\left( \log |a|\right)^{1-\rho}
\end{eqnarray*}
for some constants $C_1,C_2>0$. If the order $\rho\geq 1$ then $\left( \log |a|\right)^{1-\rho}$ is bounded above.
Consequently there are $C_3, C_4>0$ such that
$$ \Sigma ^R (\rho ,a)\geq C_3 \int_{\log R-\Delta }^{R-\Delta } \frac{T(r)}{r^{\rho +1}} dr -C_4.$$
In the case when $0<\rho<1$, we have $\left( \log |a|\right)^{1-\rho} \leq \left( \log R\right)^{1-\rho}
\asymp \left( \log (R-\Delta)\right)^{1-\rho}$. Therefore
$$ \Sigma ^R (\rho ,a)\geq C_3  \int_{\log R -\Delta }^{R-\Delta } \frac{T(r)}{r^{\rho +1}} dr -C_1
-C_5   \left( \log (R-\Delta)\right)^{1-\rho}$$
for some $C_5>0$. The assertion follows now from the assumption (\ref{eq diver}).
 \endpf

\

\ni A direct application of Theorem~\ref{theo main} gives now

\

\bcor\lab{lemm zero pressure} If $f:\C\to\oc$ is
dynamically regular,
then there exists a unique $h>\rho/\alpha$ such that $\P(h)=0$ and
$f$ has a $|f'|_\sigma^h$-conformal measure $m_h$. \ecor

\

\ni The definitions of Hausdorff measure as well as Hausdorff
dimension can be found for example in \cite{mat} or \cite{pu}. The
symbol $\H^t_\sg$ refers to the $t$-dimensional Hausdorff measure
evaluated with respect to the Riemannian metric $d\sg$. Fix
$t>\rho/\alpha$. By Theorem~\ref{theo conf meas} there exists $m_t$,
an $e^{\P(t)}|f'|_\sg^t$-conformal measure, and let $m_t^e$ be its
Euclidean version defined in the previous section. Then a
straightforward calculation shows that
\begin{equation}\lab{1030805}
{dm_t^e\circ f\over dm_t^e}(z)=e^{\P(t)}|f'(z)|^t, \  z\in J(f).
\end{equation}
Fix any radius
$$
R\in (0,\d (f)).
$$
So, if $z\in J(f)$, $n\ge 0$, and $z\in f^{-n}(w)$, then there
exists a unique holomorphic inverse branch $f_z^{-n}:D(w,4R)\to\C$
of $f^n$ sending $w$ to $z$. Recall that $K_\sg$ is the
$\sg$-adjusted Koebe constant produced in Lemma~\ref{sgbdt}. It
follows from this lemma that
\begin{equation}\lab{2030805}
D_\sg\(z,K_\sg^{-1}R|w|^{-\a_2}|(f^n)'(z)|_\sg^{-1}\) \sbt
f_z^{-n}\(D_\sg(w,R|w|^{-\a_2})\) \sbt D_\sg\(z,K_\sg
R|w|^{-\a_2}|(f^n)'(z)|_\sg^{-1}\)
\end{equation}
and that
\begin{equation}\lab{4030805}
m_t\(f_z^{-n}\(D_\sg(w,R|w|^{-\a_2})\)\) \comp
e^{-\P(t)n}|(f^n)'(z)|_\sg^{-t}m_t\(D_\sg(w,R|w|^{-\a_2})\).
\end{equation}
We recall that the radial Julia set is the set of points of $J(f)$
that do not escape to infinity:
$$
J_r(f)=\{z\in J(f): \; \liminf_{n\to\infty } |f^n(z)|<\infty\}
$$
and, obviously,
$$
J_r(f)=\bu_{M>0}J_{r,M}(f)=\bu_{M>0}\{z\in
J(f):\liminf_{n\to\infty}|f^n(z)|<M\}.
$$

\

\section{Proof of Bowen's formula}

\ni We start the proof of Bowen's formula (Theorem \ref{t1032605}) by
the following observation which, together with Lemma \ref{lemm zero
pressure}, shows in particular that $\HD (J_r(f))\leq h$.

\blem\lab{l1030805} If $t>\rho/\alpha$ such that $\P(t)\le 0$, then
$\H_\sg^t(J_r(f))<+\infty$. \elem

{\sl Proof.}
Since $\mu_t$ is an ergodic measure there is $M>0$ so large that $\mu_t(J_{r,M}(f))=1$.
Consequently $m_t(J_{r,M}(f))=1$. Since
$J(f)\cap\ov D(0,M)$ is a compact set,
$$
Q_M:=\inf\{m_t\(D_\sg(w,R|w|^{-\a_2}):w\in J(f)\cap D(0,M)\)\}>0.
$$
Now, fix $z\in J_{r,M}(f)$ and consider an arbitrary integer $n\ge
0$ such that $f^n(z)\in D(0,M)$. Recall that $D(0,T)\cap J(f)=\es$.
It follows from (\ref{2030805}) and (\ref{4030805}) that
$$
\aligned
m_t\(D_\sg\(z &,K_\sg R|f^n(z)|^{-\a_2}|(f^n)'(z)|_\sg^{-1}\)\)\gek \\
&\gek e^{-\P(t)n}|(f^n)'(z)|_\sg^{-t}m_t\(D_\sg(f^n(z),R|f^n(z)|^{-\a_2})\) \\
&\ge Q_M(K_\sg R)^{-t}e^{-\P(t)n}|f^n(z)|^{\a_2t}
    \(K_\sg R|f^n(z)|^{-\a_2}|(f^n)'(z)|_\sg^{-1}\)^t \\
&\ge Q_M(K_\sg R)^{-t}T^{\a_2t}\(K_\sg
R|f^n(z)|^{-\a_2}|(f^n)'(z)|_\sg^{-1}\)^t.
\endaligned
$$
Thus, there exists $c>0$ such that for every $z\in J_{r,M}(f)$
$$
\limsup_{r\to 0}{m_t(D_\sg(z,r))\over r^t} \ge
\limsup_{n\to\infty}{m_t\(D_\sg\(z,K_\sg
R|f^n(z)|^{-\a_2}|(f^n)'(z)|_\sg^{-1}\)\)
    \over \(K_\sg R|f^n(z)|^{-\a_2}|(f^n)'(z)|_\sg^{-1}\)^t}
\ge c.
$$
Applying now Besicovic's Covering Theorem, it immediately follows
from this inequality that $\H_\sg^t(J_{r,M}(f))\le c^{-1}$. Since
for every $x\ge M$, $m_t\(J_{r,x+1}(f)\sms J_{r,x}(f)\)=0$, an
argument similar to the one above gives that
$\H_\sg^t\(J_{r,x+1}(f)\sms J_{r,x}(f)\)=0$. Since $J_r(f)=
J_{r,M}(f)\cup \bu_{n=0}^\infty \(J_{r,M+n+1}(f)\sms
J_{r,M+n}(f)\)$, the proof is complete. \endpf

\

\ni In order to complete the proof of Bowen's formula we have to
establish that $\HD (J_r(f))\geq h$. We will do this in adapting the
corresponding proof in \cite{uz1}. The first step is to show that
$f$ has a finite and strictly positive Lyapunov exponent.

\

\blem \label{lemm lyap} We have that
$$
0< \chi
= \int \log |f'| d\mu_h
= \int \log|f'|_\sg d\mu_h
<\infty \; .
$$
\elem

{\sl Proof.} The equality $\int \log |f'| d\mu_h = \int \log
|f'|_\sigma d\mu_h$ follows from
$$
\log|f'|_\sigma(z) = \log |f'(z)|
+\alpha _2 (\log |z| -\log |f(z)|)
$$ and the $f-$invariance of
$\mu_h$. We have to prove finiteness of $\int \log f_{\sigma} '
d\mu_h$. In order to do so, consider the annulus
$A_j=D(0,2^{j+1})\setminus D(0,2^j)$. In this annulus we have
\begin{itemize}
    \item[(i)] $\mu_h (A_j) \preceq 2^{-j\delta}$ because of Lemma
    \ref{prop 2} and the fact that $d\mu_h =\psi _h dm_h$ with $\psi_h$ bounded.
    \item[(ii)] $f_{\sigma} ' (z) \preceq
    |z|^{\al}\preceq 2^{j\al}$ due to the balanced growth condition
    (\ref{eq condition beta}).
\end{itemize}
The finiteness of the integrals in the lemma follows. Finally $\chi
>0$ since $f$ is expanding.
\endpf

\

\ni We can now complete the proof of Theorem~\ref{t1032605} by
establishing the following.

\

\blem
$\HD (J_r(f))\geq h.$
\elem
{\sl Proof.} Fix $\ep >$ such that the Lyapunov exponent defined in
Lemma \ref{lemm lyap} $\chi > \ep $. Since $\mu_h(J_r(f))=1$ and
since $\mu_h$ is ergodic $f-$invariant, it follows from Birkhoff's
ergodic theorem and Jegorov's theorem that there exists a Borel set
$Y\subset J_r(f)$ and an integer $K\geq 1$ such that $\mu_h(Y)\geq
{1\over 2}$ and such that for every $z\in Y$ and $n\geq k$
\begin{equation} \label{eq lyap1}
\left| {1\over n} \log |(f^n) ' (z)|  -\chi \right| <\ep \;\; and \;\;
\left| {1\over n} \log |(f^n) ' (z)|_\sigma \  -\chi \right| <\ep.
\end{equation}
Let $R=\dist(\post , J(f))/4$. Given $z\in Y$ and $r\in (0,R)$, let
$n\geq 0$ be the largest integer such that
$$ D(z,r) \subset f_z^{-n} (D(f^n(z),R)).$$
There is $r_z>0$ such that for any $0<r<r_z$ the integer $n$ defined
above is $n\geq k$. By the definition of $n$, $D(z,r)$ is not
contained in $f_z^{-(n+1)} (D(f^{n+1}(z),R))$. Koebe's distortion
theorem yields now
\begin{equation}\label{eq 11} r \leq KR
|(f^n)'(z)|^{-1} \;\;\; and \;\; \;r  \geq K^{-1} R
|(f^{n+1})'(z)|^{-1}.
\end{equation}
Passing to the $h$-conformal measure $m_h$, we get from
(\ref{eq lyap1}) that
$$
\aligned
m_h(D(z,r))&\leq m_h\(f_z^{-n} (D(f^n(z),R)\)
\asymp |(f^n)'(z)|_\sigma ^{-h} m_h \(D(f^n(z),\delta)\)\\
&\leq |(f^n)'(z)|_\sigma^{-h} \leq e^{-hn (\chi -\ep )}.
\endaligned
$$
On the other hand, (\ref{eq lyap1}) together with (\ref{eq 11})
give
$$e^{-(n+1) (\chi +\ep )} \leq |(f^{n+1})'(z)|^{-1} \preceq r.$$
Therefore
$$m_h(D(z,r)) \preceq r^{h\(\frac{n+1}{n}\frac{\chi -\ep}{\chi+\ep}\)}.$$
When $r\to 0$ then $n=n(r)\to \infty $ from which we get that
$$\limsup_{r\to 0} \frac{m_h(D(z,r))}{r^{h-\ep'}}\preceq 1$$
for every $\ep '>0$. This gives $\HD (J_r(f))\geq
h-\ep'$ and the lemma follows in taking $\ep' \to
0$.
\endpf


\section{Real analyticity of the hyperbolic dimension} \label{sec 8}

\ni In this section we prove Theorem~\ref{1.1}. From now on we suppose $\a_1\geq0$.

\subsection{J-stability} \label{sec J-stability}
 The work of Lyubich and Ma\~{n}\'e-Sad-Sullivan \cite{l1, mss}
on the structural stability of rational maps has been generalized to
entire functions of the Speiser class by Eremenko-Lyubich \cite{el}.
Note also that they show that any entire function of the Speiser
class is naturally imbedded in a holomorphic family of functions in
which the singular points are local parameters.

Here we collect and adapt to the meromorphic setting the
facts that are important for our needs. We also deduce from the bounded deformation assumption of $\cM_\La$ near
$\fo$ a bounded speed condition of the involved holomorphic
motions. A \it holomorphic motion \rm of a set $A\subset \C$ over
$U$ originating at $\l^0$ is a map $h : U \times A \to \C$
satisfying the following conditions:
\begin{enumerate}
    \item The map $\l \mapsto h(\l , z) $ is holomorphic for every
    $z\in A$.
    \item The map $h_\l : z\mapsto h_\l (z) = h(\l ,z) $ is
    injective for every $\l \in U$.
    \item $h_{\l^0} =id$.
\end{enumerate}
The $\l$--lemma  \cite{mss} asserts that such a holomorphic motion
extends in a quasiconformal way to the closure of $A$. Further
improvements, resulting in the final version of Slodkowski
\cite{sk}, show that each map $h_\l$ is the restriction of a global
quasiconformal map of the sphere $\oc$.
\sp\ni Let us call $\fo\in \cM_\La$ \it holomorphically J-stable \rm
 if there is a neighborhood $U\subset \La$ of $\l^0$ and a holomorphic motion
 $h_\l$ of $\julo$ over $U$ such that $h_\l (\julo ) = {\mathcal J}(f_\l )$ and
$$
h_\l \circ \fo = f_\l \circ h_\l \quad  on  \;\;\julo
$$
for every $\l \in U$.

\

\blem \label{3.1} A function $\fo\in \cM_\La$ is holomorphically
J-stable if and only if, for every singular value $a_{j,\l^0} \in
sing(\fo^{-1})$, the family of functions
$$\l \mapsto f_\l^n (a_{j,\l })\;, \;\; n\geq 1,$$
is normal in a neighborhood of $\l^0$.
\elem

\bpf This can be proved precisely like for rational functions
because the functions in the Speiser class $\cS$ do not have wandering
nor Baker domains
(see \cite{l2} or \cite[p. 102]{bm}). \epf

 From this criterion together with the description of the
components of the Fatou set one easily deduces the following.

\blem \label{3.2} Each $\fo\in \cH\cM_\La$ is holomorphically $J$--stable
and $\cH\cM_\La$ is open in $\cM_\La$. \elem

 We now investigate the speed of the associated holomorphic
motion.

\bprop \label{3.3} Let $\fo\in \cH\cM_\La$ and let $h_\l$ be the associated
holomorphic motion over $U\subset \Lambda$ (cf. Lemma \ref{3.2}). If
$\cM_U$ is of bounded deformation, then there is $C>0$ such that
$$
\left|\frac{\partial h_\l (z)}{\partial \l _j} \right| \leq C
$$
for every $z\in \julo$ and $j=1,...,N$. It follows that $h_\l$
converges to the identity map uniformly on $\julo $ and,
replacing $U$ by a smaller neighborhood if necessary, that there
exists $0<\tau \leq 1$ such that $h_\l$ is $\tau$-H\"older for every
$\l \in U$. \eprop

\bpf Let $h_\l$ be the holomorphic motion such that $f_\l \circ h_\l
= h_\l \circ \fo$ on $\julo$ for $\l \in U$ and such that there are
$c>0$ and $\rho >1$ for which
\begin{equation} \label{3.4}  |(f_\l^n)'(z)| \geq c\rho ^n \quad for
\; every\; n\geq 1, \; z\in \cJ_{f_\l} \; and \; \l \in
U.\end{equation}
 (cf. Fact \ref{2.2}; this is the only place where $\a_1\geq 0$ is used). Denote $z_\l =h_\l (z)$ and consider
$$ F_n (\l , z) = f_\l ^n (z_\l ) - z_\l.$$
The derivative of this function with respect to $\l_j$ gives
$$ \frac{\partial}{\partial \l_j}F_n(\l,z)
=\frac{\partial f_\l^n}{\partial \l_j}(h_\l(z)) + (f_\l^n)'
(h_\l(z)) \frac{\partial}{\partial \l_j}h_\l(z)
-\frac{\partial}{\partial\l_j} h_\l(z).$$ Suppose that $z$ is a
repelling periodic point of period $n$. Then $\l \mapsto
F_n(\l,z)\equiv 0$ and it follows from (\ref{3.4}) that
$$\left|\frac{\partial h_\l(z)}{\partial \l_j}\right| = \left|
\frac{\frac{\partial f_\l^n}{\partial
\l_j}(z_\l)}{1-(f_\l^n)'(z_\l)} \right| \preceq \left|
\frac{\frac{\partial f_\l^n}{\partial \l_j}(z_\l)}{(f_\l^n)'(z_\l)}
\right| = \De _{n,j}.$$ Since $\frac{\partial f_\l^n}{\partial
\l_j}(z_\l)= \frac{\partial f_\l}{\partial
\l_j}\big(f_\l^{n-1}(z_\l)\big) +f_\l'\big(f_\l^{n-1}(z_\l)\big)
\frac{\partial f_\l^{n-1}}{\partial \l_j}(z_\l)$ we have
$$\De_{n,j} \leq \frac{\left|\frac{\partial f_\l}{\partial
\l_j}(f^{n-1}_\l(z_\l))\right|}{\left|f_\l'(f^{n-1}_\l(z_\l))\right|}
\frac{1}{\left|(f_\l^{n-1})'(z_\l)\right|} + \De_{n-1,j}.$$ Making
use of the expanding (\ref{3.4}) and the bounded deformation
(\ref{1.3}) properties it follows that
$$\De_{n,j} \leq \frac{M}{c\rho^{n-1}} + \De_{n-1,j}.$$
The conclusion comes now from the density of the repelling cycles in
the Julia set $\julo$:
$$ \left| \frac{\partial h_\l(z)}{\partial \l_j}\right|\preceq
\frac{M}{c}\frac{\rho}{\rho-1} \;\; for \;\; every \; z\in \julo.$$
The H\"older continuity property is now standard (see \cite{uz1}).
\epf

\

Concerning the divergence type condition and the growth condition on
the characteristic function (\ref{eq diver}), these are stable in the 
sense that if $f_{\l _0}$ satisfies it, then $f_\l$ also has this property for all
$\l$ in some neighbourhood of $\l_0$. For example, if $f_{\l _0}$ is entire,
then 
$$ T(r, {f_{\l _0}})= m(r,\infty ) = \frac{1}{2\pi} \int _{-\pi}^\pi \log ^+ |f_{\l _0}(re^{i\theta})|\, d\theta.$$
If (\ref{eq diver}) holds for this function, then it follows from that
expression and from the uniform convergence on compact sets of $f_\l$ to ${f_{\l _0}}$
that (\ref{eq diver}) does hold for all $f_\l$, $\l$ in some neighborhood of $\l_0$.

\

\subsection{The spectral gap of the (real) transfer operator.}

In order to get the necessary spectral properties of the transfer
operator, one does work with the space of H\"older continuous
functions $\cH_\tau = \cH_\tau (\jul , \C )$, $0<\tau \leq 1$.
However, the function $|f'|_\sigma ^{-1}$ is not necessary in this space.
It follows from the distortion property (\ref{4.3}) below that it belongs to the
following slightly more general one. In order to introduce it consider
$w\in \jul$ and denote the $\tau $--variation of a function
$g:\jul \cap D(w,\d ) \to \C$ by
$$v_{\tau , w}(g)  = \sup \left\{ \frac{|g(x) -g(y)| }
{|x-y|^\tau} \; ;  \;\;  x,y\in \jul \cap D(w,\d )
\right\}.$$
The H\"older space $H_\tau$ we work with consists in bounded functions
$g:\jul \to \C$ such that $v_{\tau , w}(g\circ f_a^{-1})$ is bounded uniformly
in $w\in \jul$ and $a\in f^{-1}(w)$. Denote
$$v_\tau (g) = \sup_{w\in \jul} \sup_{a\in f^{-1}(w)}v_{\tau , w}(g\circ f_a^{-1}).$$
 The space $H_\tau$ endowed with the
norm $\|g\|_\tau = v_\tau (g) +\|g\|_\infty$ is a Banach space
densely contained in $C_b$. Here is the classical estimation which
is based on the H\"older property of $g\in H_\tau $ and the
expanding property of $f$:
$$\aligned
\left|  \npft^n g(z) - \npft^n g(w) \right|& = e^{-nP(t)}\Big|\sum
_{a\in f^{-n}(z)} |(f^n)'(a)|_\sg ^{-t}g (a)  - \sum _{b\in
f^{-n}(w)}
|(f^n)'(b)|_\sg ^{-t}g (b)\Big|\\
&\leq  e^{-nP(t)} \sum _{a\in f^{-n}(z)} |(f^n)'(a)|_\sg
^{-t}\left|g(f^{-n}_a(z)) -g(f^{-n}_a(w))\right|\\
&+e^{-nP(t)} \sum _{b\in f^{-n}(w)}
 \big| |(f^n)'(f^{-n}_b(z))|_\sg ^{-t}
  -|(f^n)'(b)|_\sg ^{-t}
  \big|\big|g(b) \big|\\
  &= I + II
\endaligned
$$
for $z,w\in \jul$ with $|z-w|<\d =\d (f)$ and where $f^{-n}_a$ is
the inverse branch of $f^{-n}$ defined on $D (z,\d )$ such that
$f^{-n}_a(z) = a$. The majorization of the first term goes as
follows:
$$\aligned
I&\leq v_\tau (g) e^{-nP(t)} \sum _{a\in f^{-n}(z)} |(f^n)'(a)|_\sg
^{-t}\left|f^{-(n-1)}_{f(a)}(z) -f^{-(n-1)}_{f(a)}(w)\right|^\tau \\
&\leq v_\tau (g) \| \npft ^n \| _\infty \sup_{a\in f^{-n}(z)}
|f^{-(n-1)}_{f(a)}(z)-f^{-(n-1)}_{f(a)}(w)|^\tau \preceq v_\tau (g) \rho
^{-(n-1)\tau}|z-w|^\tau .
\endaligned
$$
Concerning the second part, one has to observe that Koebe's
distortion theorem implies that for any $n\geq 1$ and $z,w,\in \jul$
with $|z-w|<\d (f)$ \begin{equation} \label{4.3}
 \big| |(f^n)'(f^{-n}_a(z))|_\sg
^{-t} -|(f^n)'(f^{-n}_a(w))|_\sg ^{-t} \big| \preceq
|(f^n)'(f^{-n}_a(z))|_\sg ^{-t} |z-w| \end{equation} where $a\in
f^{-n}(z)$. Therefore
$$II\preceq e^{-nP(t)}\sum_{b\in f^{-n}(w)} |(f^n)'(b)|_\sigma ^{-t} |z-w| |g(b)| \preceq
 \|\npft ^n\| \|g\|_\infty |z-w|.$$
Altogether we have \begin{equation} \label{4.4} \left|  \npft^n g(z)
- \npft^n g(w) \right|\preceq \big(\rho ^{-(n-1)\tau } v_\tau(g)
+\|g\|_\infty \big) |z-w|^\tau\end{equation} for all $z,w\in \jul$
with $|z-w|<\d (f)$. We proved

\

\blem \label{4.5} $ \npft (H_\tau ) \subset H_\tau$ and, for any
$g\in H_\tau$ and $n\geq 1$,
$$ \|\npft ^n (g)\|_\tau\preceq
\rho^{-(n-1)\tau } v_\tau (g) + \|g\|_\infty.$$ \elem

 \

 \ni If $B $ is a
bounded subset of $H_\tau$ then (\ref{4.4}) and the fact that
$\|\npft^n\|_\infty$ is uniformly bounded yields that $\cF=\{\npft
(g); \; g\in B \}$ is a equicontinuous bounded subfamily of $(C_b,
\| .\|_\infty )$ . The following observation follows then precisely
like in \cite[Lemma 4.2]{uz1} (using $\lim_{w\to\infty } \pft \1 (w)
=0$ which is Lemma \ref{lemm 2}).

\

\blem \label{4.6} If $B$ is a bounded subset of $H_\tau$, then
$\npft (B) $ is a precompact subset of $(C_b, \| .\| _\infty)$.
\elem

\

\ni We are now in the position to apply Ionescu-Tulcea and
Marinescu's Theorem 1.5 in \cite{im}. Combined with \cite{du2} (see
\cite{uz1} were these facts are explained in detail) we finally get:
\ \bprop \label{4.7} For all $t>\rho /\a $ there is $r\in (0,1)$
such that the spectrum $\sigma (\npft ) \subset \D (0,r)\cup \{1\}$
and the number $1$ is a simple isolated eigenvalue of the operator
$\npft$ of $H_\tau$. \eprop

\

\subsection{Complexified transfer operator} In the remainder of the
paper we consider a hyperbolic function $\fo \in \cH\cM_\La$. Let $U\subset
\Lambda $ be a neighborhood of $\l _0$ on which $f_\l$ is hyperbolic
and holomorphically $J$--stable and let
$$
\pftl g(w) = \sum_{z\in f_\l ^{-1} (w)} |f_\l ' (z)|_\sg ^{-t} g(z)
\; , \;\; t> \frac{\rho }{\a} ,
$$
be the induced family of (real) transfer operators acting
continuously on $C_b (\jull , \C )$ and on $H_1 (\jull , \C )$. In
order to be able to work on the fixed Julia set $\julo$ we conjugate
these operators by $T_\l : C_b (\jull , \C ) \to C_b (\julo , \C )$
where $T_\l (g) =g\circ h_\l $ and where $h_\l$ is the associated
holomorphic motion. Put
$$
L(t,\l )= T_\l \circ \pftl \circ T_\l ^{-1}
$$
to be the resulting bounded operator of $C_b = C_b (\julo,\C )$. We
have that
$$ L(t,\l ) (g) (w)
=\sum_{z\in \fo ^{-1} (w)} |f_\l ' (h_\l (z))|_\sg ^{-t} g(z) \; ,
\;\; w\in \julo \; , \;\; g\in C_b.
$$
Our aim is to establish real analyticity of the hyperbolic dimension
of $f_\l$. In order to do so we have to embed these operators in a
holomorphic family
$$
(t,\l ) \in \C \times \C^{2d} \rightarrow L(t,\l ) \in L(H_\tau).
$$
In order to do so, we follow \cite{uz1} and start with complexifying
the potentials $|f_\l '|_\sg ^{-t} \circ h_\l $. Denote again
$z_\l=h_\l (z)$, $ z\in \julo$ and $\l\in \D_{\C^d}(\l^0, R)$. Remember that $h_\l \to id$
uniformly in $\julo$ (Proposition \ref{3.3}). Since $0\notin\julo$ the function
$$\Psi_z(\l ) = \frac{f_\l '(z_{\l})}{f_{\l^0}'(z)}\lt(\frac{z_\l}{z}\rt)^{\a_2}\lt(\frac{f_{\l^0}(z)}{f_\l (z_\l)}\rt)^{\a_2}$$
is well defined on the simply connected domain $\D_{\C^d}(\l^0, R)$.
Here we choose $w\mapsto w^{\a_2}$ so that this map fixes $1$ which implies that
$$\Psi_z (\l^0)=1 \quad for \; every \;\; z\in \jo =  \julo\sms \fo ^{-1}(\infty).$$
For this function one has the following uniform estimate.

\blem\label{5.1} For every $\e >0$ there is $0<r_\e <R$ such that
$|\Psi _z (\l ) -1| < \e$ for every $\l\in \D_{\C^d}(\l^0,r_\e)$ and every
$z\in \jo$. \elem

\bpf
Suppose to the contrary that there is $\ep>0$ such that for some
$r_j\to 0$ there exists $\l_j\in \D_{\C^d}(\l^0,r_j)$ and $z_j\in \jo$
with $|\Psi_{z_j}(\l_j) -1|>\ep$. Then the family of functions
$$\cF =\{\Psi_z \, ; z\in \jo \}$$
cannot be normal on any domain $\D_{\C^d}(\l^0 ,r)$, $0<r<R$. This is however not true.
Indeed, the
balanced growth condition (\ref{eq
condition beta}) yields
$$
|\Psi_z(\l) | \leq \kappa^2\lt| \frac{z_\l}{z} \rt|^\a \quad for \; every \; z\in \jo \; and \; |\l-\l^{0}|<R.
$$
Since $h_\l \to Id$ uniformly in $\C$ it follows immediately that $\cF$ is normal on some
disk $\D_{\C^d}(\l^0 ,r)$, $0<r<R$.
 \epf

\

\ni We can now proceed precisely as in \cite{uz1} (or in
\cite{cs2}). Embed $$\l =(\l_{d-1} , ..., \l _0 ) = (x_{d-1}
+iy_{d-1},..., x_0+iy_0) \in \C^d$$ into $\C^{2d}$ by the formula
$\l \mapsto (x_{d-1}, y_{d-1},..., x_0, y_0) \in \C^{2d}$, replace
in the power series of the, for every $z\in\jo$, real analytic
functions
$$(t,\l ) \mapsto |f'_{\l } (z_\l )|_\sg ^{-t} =
\exp\{-t \Re \log f'_{\l , \sg} (z_\l )\} = |\fo' (z )|_\sg ^{-t}
\exp\{-t \Re \log \Psi _z (\l )\},$$
 $\| \l -\l _0\|<R$ and $\Re (t)> \frac{\rho }{\a }$, the real
 numbers $x_j = \Re \l _j$, $ y_j = \Im \l _j$ by complex numbers
 and obtain by a straightforward adaption of the arguments given in
 \cite{uz1, cs2} the following:

\

\bprop \label{5.2} There is $R>0$ such that, for every $z\in \jo$,
the function
$$(t,\l ) \mapsto  \ph _{t,\l } (z) =|\fo' (z )|_\sg ^{-t}
\exp\{-t \Re \log \Psi _z (\l )\}$$ can be extended to a holomorphic
function on $\{\Re t>\frac{\rho}{\a }\}\times \D _{\C ^{2d}} (\l _0
, R)$. In addition, this extension that we still denote $\ph _{t,\l
}$ has the following properties:
\begin{enumerate}
    \item $|\ph _{t,\l } (z) | \asymp |\fo' (z )|_\sg ^{-t}$.
    \item There is $0<\tau \leq 1$ such that $\ph _{t,\l } \in
    H_\tau $ and $(t,\l ) \mapsto \ph _{t,\l } \in H_\tau$ is
    continuous.
    \item $\ph _{t,\l }$ is uniformly dynamically H\"older.
\end{enumerate}
 \eprop

\

\ni A continuous function $\ph : \julo \to \C$ is called
$c_\ph$--\it dynamically H\"older \rm of exponent $\tau$ if
$$
|\ph_n ((\fo)^{-n}_a(z)) - \ph_n((\fo)^{-n}_a(w)) | \leq c_\ph
|\ph_n ((\fo)^{-n}_a(z)) ||z-w|^\tau
$$
for $a\in \fo^{-n}(z)$, $|z-w|<\d (\fo )$ and with $\ph_n (a) = \ph
(a) \ph(\fo (a)) \cdot...\cdot\ph (\fo^{n-1} (a))$. As we noted in
(\ref{4.3}), $\ph_{t,\l^0 }(z) = |\fo' (z )|_\sg ^{-t}$ is
dynamically H\"older. The family of potentials $\ph_{t,\l }$ is
called \it uniformly dynamically H\"older \rm if the involved
constants $\tau , c_\ph $ above can be chosen to be valid for all
the potentials of the family. Item (1) of the preceding proposition
means in particular that the transfer operators \begin{equation}
\label{5.3} L(t,\l ) (g)(w) = \sum_{z\in \fo ^{-1}(w)} \ph_{t,\l }
(z) g(z) \end{equation}
 are (uniformly) bounded on $C_b$ (such
potentials are also called (uniformly) summable). In fact, much more
is true since Proposition \ref{5.2} together with Corollary 7.7 of
\cite{uz1} yield:

\

\bcor \label{5.4} There are $0<\tau \leq 1$ and $R>0$ such that the
operators $L(t, \l )$ are bounded operators of $H_\tau $ and such
that the map
$$ (t,\l ) \in \{\Re t >\rho/\a\}\times \D _{\C^{2d}} (\l
_0 ,R) \mapsto L(t,\l ) \in L(H_\tau )$$ is holomorphic. \ecor

\

\subsection{Real analyticity of the hyperbolic dimension} We are now
in position to proof Theorem~\ref{1.1}. We take the notation of the
preceding section, in particular $\fo\in \cH$ is a hyperbolic
function. Consider a real $t_0>\frac{\rho }{\a }$. Then we have
$\pftolo = L(t_0,\l _0)\in L(H_\tau )$ and this operator has a
simple and isolated eigenvalue which is $\g (t_0 ,\l _0) =
e^{P_{\l^0}(t_0)}$, where $P_{\l^0}(t_0)$ is the topological
pressure of $\fo$ at $t_0$ (see Proposition~\ref{4.7}). From the
perturbation theory for linear operators (see \cite{ka}) it follows
now that there is $r>0$ and a holomorphic map
$$
(t,\l )\in \D_\C (t_0 ,r)\times \D_{\C^{2d}}(\l^0 , r) \mapsto \g
(t, \l)
$$
such that
\begin{enumerate}
    \item $\g (t,\l )$ is a simple isolated eigenvalue of $L(t, \l )\in
L(H_\tau )$ and
    \item there is $\b>0$ such that the spectrum $$\sg (L(t,\l ))\cap
    \D (e^{P_{\l^0}(t_0)},\b ) = \{\g (t,\l )\}$$ for all
    $(t,\l )\in \D_\C (t_0  , r)\times \D_{\C^{2d}}(\l^0 ,
    r)$.
\end{enumerate}
Coming now back to the initial parameters, real $t$ and $\l \in
\C^d$, we remember that the operators $L(t,\l )$ are conjugate to
$\pftl$ via the operator $T_\l$ that consist in composition with the
holomorphic motion $h_\l$. From the H\"older continuity property
(Proposition \ref{3.3}) of $h_\l$ it follows that we may assume that
there is $0<\tau \leq 1$ such that $T_\l (H_1(\jull , \C) ) \subset
H_\tau (\julo , \C)$ for all $\|\l - \l _0\| <r$. Consequently
$e^{P_\l (t)}$, $P_\l (t)$ the topological pressure of $f_\l$ at
$t$, is an eigenvalue of $\pftl$ provided we can show the following:

\

\blem \label{6.1} For every $t>\rho/\a$ the function $\l
\mapsto \P_\l (t)$ is continuous \elem

\bpf We have that $\P_\l (t)=\lim_{n\to\infty } \frac{1}{n}\log
\sum_{z\in f_\l ^{-1}(w)} |(f_\l ^n)'(z)|_\sg ^{-t}$ with $w\in \jull
$ is any finite point. The continuity assertion results directly
from Lemma \ref{5.1} since it is shown there that for any $z\in \jo$
$$ (1-\e )^n
\leq\frac{|(f_\l ^n)'(h_\l (z))|_\sg}{ |(\fo ^n)'(z)|_\sg}
\leq (1+\e)^n.
$$
\epf

 \

 \ni Altogether we obtained real analyticity of the pressure
function. From Bowen's formula (Theorem \ref{t1032605})) we know
that the hyperbolic dimension $\HD (\radl )$ is the only zero of the
pressure function $t\mapsto \P_\l (t)$. Real analyticity of this zero
with respect to $\l$ results from the implicit function theorem
since clearly
$$\frac{\partial}{\partial t} \P_\l (t) \leq -\log \rho <0$$
where $\rho >1$ is the expanding constant that is common to the
$f_\l$ (see Fact \ref{2.2}).

\section{Around Theorem~\ref{1.4}\label{proofoft1.1}}
\ni In this section we derive the most transparent consequences of
Theorem~\ref{1.1}, notably Theorem~\ref{1.4}. We begin with the
following.

\

\bthm\label{t1072005}
Let $f_\l = f\circ P_\l$ with $f:\C\to\oc$ a
meromorphic function and, for every $\l=(\l_d,\l_{d-1},\ld,\l_1,\l_0)
\in\C^{d+1}$, $P_\l:\C\to\C$ is the polynomial given by the formula
$P_\l(z)=\sum_{j=0}^d\l_jz^j$. Suppose that $f_{\l^0}$ is hyperbolic
and that there is a neighborhood $U\subset \C\sms \{0\}\times \C^d$
of $\l^0$ such that $\{f_\l ; \; \l \in U\}$ is uniformly
balanced with $\a_2\ge 1$ and $\a_1\geq 0$. Then the function $\l\mapsto
\HD(\JPl)$ is real-analytic near $\l^0$.
\ethm
\bpf Put $f_\l=f\circ P_\l$. For every $\g\in\C^{d+1}$ put
$$
Q_\g(z)=z^d+\sum_{j=0}^{d-1}\g_j\g_d^{-j}z^j  \ \text{ and }  \
g_\g=\g_df\circ Q_\g.
$$
Consider also $H$, the change of coordinates in the parameter space,
given by the formula
$$
H(\l_d,\l_{d-1},\ld,\l_1,\l_0)=(\l_d^{1/d},\l_{d-1},\ld,\l_1,\l_0),
$$
where $\l_d\mapsto\l_d^{1/d}$ is a holomorphic branch of $d$th
radical defined on the ball $\D_{\C^d}(\l^0_d,|\l^0_d|)$. Let $T_\g: \C\to\C$
be the multiplication map defined as $T_\g(z)=\g_d^{-1}z$. Notice
that
\begin{equation}\label{1072005}
T_{H(\l)}\circ g_{H(\l)}\circ T^{-1}_{H(\l)}=f_\l.
\end{equation}
So, $\Jl=T_{H(\l)}({\mathcal J}_r(g_{(H(\l)}))$, and in consequence,
$$
\HD(\Jl)=\HD({\mathcal J}_r(g_{(H(\l)})).
$$
Since in addition
$H(\l^0)=\((\l_d^0)^{1/d},\l^0_{d-1},\ld,\l^0_1,\l^0_0\)$, in order
to prove our theorem, it is enough to show that the map $\g\mapsto
\HD({\mathcal J}_r(g_\g))$ is real-analytic near the point $\g^0=
\((\l_d^0)^{1/d},\l^0_{d-1},\ld,\l^0_1,\l^0_0\)\in H(U)$. It follows from
(\ref{1072005}) that for every $\l$ in a neighbourhood of $\l^0$ and
every $z\in {\mathcal J}(g_{(H(\l)}))$, we have
$$
|g_{H(\l)}(z)| =|\l_d^{1/d}||f_\l(T_{H(\l)}(z))|  \
 \text{ and }  \
|g_{H(\l)}'(z)|=|f_{\l}'(T_{H(\l)}(z))|.
$$
Consequently,
$$
{|g_{H(\l)}'(z)|\over|g_{H(\l)}(z)|} =|\l_d^{-{1\over
d}}|{|f_{\l}'(T_{H(\l)}(z))|\over|f_{\l}(T_{H(\l)}(z))|}.
$$
Since $T_{H(\l)}(z)\in\jull$ and since $|\l_d^{-{1\over d}}|$ is
bounded away from zero and infinity on a neighbourhood of $\l^0$, it
follows from the uniform balanced growth of $\{f_\l ; \, \l \in U\}$
that $g_{\l}$ also has this property for $\l $ near $\l^0$. Aiming
to apply Theorem~\ref{1.1}, we are therefore left to show that for a
sufficiently small bounded neighbourhood of $\g^0$, the family
$\cM_U=\{g_\g\}_{\g\in H(U)}$ is of bounded deformation. We have for
every $z\in\C$ that
\begin{equation}\label{2072005} g_\g'(z) =\g_df'(Q_\g(z))Q_\g'(z)
=\g_df'(Q_\g(z))\lt(dz^{d-1}+\sum_{j=1}^{d-1}j\g_j\g_d^{-1}z^{j-1}\rt),
\end{equation} \begin{equation}\label{4072105} \aligned \dot g_d(\g,z) &:={\bd
g_\g\over\bd\g_d}(z)
=f(Q_\g(z))+\g_df'(Q_\g(z)){\bd Q_\g(z)\over\bd\g_d}, \\
&=f(Q_\g(z))+\g_df'(Q_\g(z))\sum_{j=1}^{d-1}-j\g_j\g_d^{-j-1}z^j,
\endaligned
\end{equation}
and
\begin{equation}\label{2072105}
\dot g_i(\g,z)
:={\bd g_\g\over\bd\g_i}(z)
=\g_df'(Q_\g(z)){\bd Q_\g(z) \over\bd\g_i}
=\g_df'(Q_\g(z))\g_d^{-i}z^i
\end{equation}
for all
$i=0,1,\ld,d-1$. Taking $U$ sufficiently small, there clearly exists
$p\in(0,+\infty)$ such that
\begin{equation}\label{1072105}
|Q_\g'(z)|\ge 1
\end{equation}
for all $\g\in H(U)$ and all $z\in\C$
with $|z|\ge p$. Now, since $g_\g$ is of uniformly rapid derivative
growth on $H(U)$ and since, after a conjugation by translation,
there exists $R>0$ such that
\begin{equation}\label{5072105}
\Jgg\cap D(0,R)=\es
\end{equation}
for all
$\g\in H(U)$, it follows from (\ref{2072005}) that $Q_\g'(z)\ne 0$ for
all $\g\in H(U)$ and all $z\in \Jgg$. By a standard compactness
argument, it then follows from $J$-stability of $g_{\g_0}$ that
decreasing $U$ appropriately, we get
$$
A:=\inf\{|Q_\g'(z)|:\g\in H(U),\, z\in \ov D(0,p)\cap\Jgg\}>0.
$$
Combining this and (\ref{1072105}), we obtain
$$
B:=\inf\{|Q_\g'(z)|:\g\in H(U),\, z\in\Jgg\}\ge \min\{1,A\}>0.
$$
It follows from (\ref{2072005}) and (\ref{2072105}) that for all
$i=0,1,\ld,d-1$ we have
\begin{equation}\label{3072105}
{|\dot g_i(\g,z)|\over |g_\g'(z)|}
=|\g_d|^{-i}{|z|^i\over|Q_\g'(z)|}
=|\g_d|^{-i}{|z|^i\over
  \lt|dz^{d-1}+\sum_{j=1}^{d-1}j\g_j\g_d^{-1}z^{j-1}\rt|}.
\end{equation}
Since obviously, $\lim_{z\to\infty}(|z|^i/|Q_\g'(z)|)\le 1/d$
uniformly with respect to $\g\in H(U)$ for all
$i=0,1,\ld,d-1$, invoking the definition of $B$, we see that
$$
B_1:=\max_{0\le i\le d-1}\lt\{
     \sup\lt\{{|z|^i\over |Q_\g'(z)|}:\g\in H(U),\, z\in\Jgg\rt\}\rt\}<+\infty.
$$
Combining this and (\ref{3072105}), we see that with $U$
sufficiently small,
\begin{equation}\label{8072105}
B_2:=\max_{0\le i\le d-1}\lt\{\sup\lt\{
     {|\dot g_i(\g,z)|\over |g_\g'(z)|}:\g\in H(U),\, z\in\Jgg\rt\}\rt\}<+\infty.
\end{equation}
It follows from (\ref{2072005}) and (\ref{4072105}) that
\begin{equation}\label{7072105}
\aligned {|\dot g_d(\g,z)|\over |g_\g'(z)|}
&:=\lt|{f\circ Q_\g(z)\over\g_d(f\circ Q_\g)'(z)}+
  {{\bd Q_\g(z)\over \bd\g_d}\over |Q_\g'(z)|}\rt| \\
&\le |\g_d|^{-1}{|f\circ Q_\g(z)|\over|(f\circ Q_\g)'(z)|}
  +\lt|{\sum_{j=1}^{d-1}-j\g_j\g_d^{-j-1}z^j\over dz^{d-1}+
   \sum_{j=1}^{d-1}j\g_j\g_d^{-1}z^{j-1}}\rt|.
\endaligned
\end{equation}
Since we have the uniformly balanced growth property, since
$\a_1\ge 0$ and $\a_2\ge 1$, and taking into account (\ref{5072105}),
we conclude that \begin{equation}\label{6072105}
B_3:=\sup\lt\{|\g_d|^{-1}{|f\circ Q_\g(z)|\over|(f\circ Q_\g)'(z)|}
   :\g\in H(U),\, z\in\Jgg\rt\}<+\infty.
\end{equation}
Since obviously,
$$
\lim_{z\to\infty}\sup_{\g\in H(U)}\lt\{\lt|
     {\sum_{j=1}^{d-1}-j\g_j\g_d^{-j-1}z^j\over dz^{d-1}+
     \sum_{j=1}^{d-1}j\g_j\g_d^{-1}z^{j-1}}\rt|\rt\}
     <+\infty,
$$
invoking the definition of $B$, we see that
$$
B_4:=\sup\lt\{\lt|
 {\sum_{j=1}^{d-1}-j\g_j\g_d^{-j-1}z^j\over dz^{d-1}+
  \sum_{j=1}^{d-1}j\g_j\g_d^{-1}z^{j-1}}\rt|:\g\in H(U),\, z\in\Jgg\rt\}
  <+\infty.
$$
Combining this, (\ref{6072105}), (\ref{7072105}), and
(\ref{8072105}), we see that
$$
\max_{0\le i\le d}\lt\{\sup\lt\{
     {|\dot g_i(\g,z)|\over |g_\g'(z)|}:\g\in H(U),\, z\in\Jgg\rt\}\rt\}<+\infty.
$$
We are done. \epf

\

\ni Note that if $d=1$, then with the notation of the proof of the
previous theorem, $g_{\l_d}=\l_df$ and, as an immediate consequence
of this proof, we have the following.

 \bcor\label{c1072105}
Suppose that $f:\C\to\oc$ is a meromorphic function and consider the
analytic family $\F=\{\l f\}_{\l\in\C\sms\{0\}}$. If $f\in\F$ is
hyperbolic and if this family is uniformly balanced near
$f$ with $\a_2\geq 1$ and $\a_1\geq 0$, then the function $\l\mapsto\HD({{\mathcal J}_r(\l f)})$ is
real-analytic in a neighbourhood of $\l^0=1$. \ecor

\

\brem\label{r2072105} If in the formulation of
Theorem~\ref{t1072005} the parameter $\l_d$ is kept fixed equal to
$1$, then the derivative $\dot g_d(\g,z)$ disappears and it suffices
to assume that $\a_2>0$ (and $\a_1\ge 0$). \erem

 We end this section by noting that Theorem~\ref{1.4} is an immediate
consequence of Theorem~\ref{t1072005} and Lemma~\ref{lemm unif balanced}.



\begin{thebibliography}{MMMM}

\bibitem[Ba]{bar} K. Bara\'nski, Hausdorff dimension and measures on Julia sets of some
meromorphic functions, Fund. Math. 147 (1995), 239-260.
\bibitem[Bw]{b} W. Bergweiler, \em Iteration of meromorphic functions, \rm Bull. A.M.S. 29:2 (1993), 151-188.
\bibitem[BM]{bm} F. Berteloot, V. Mayer, \em Rudiments de dynamique holomorphe,   \rm
Cours sp\'ecialis\'es 7, SMF (2000).

\bibitem[Br]{borel} \'E. Borel, \em Sur les z\'eros des fonctions enti\`eres,   \rm  Acta Math. 20 (1897), 357-396.
\bibitem[Bw]{bow} R. Bowen, \em Hausdorff dimension of quasi-circles, \rm Publ. Math. IHES, 50
(1980), 11-25.
\bibitem[CS1]{cs1} I. Coiculescu, B. Skorulski, \em Thermodynamic formalism of transcendental
entire maps of finite singular type, \rm Preprint 2004.
\bibitem[CS2]{cs2} I. Coiculescu, B. Skorulski, \em Perturbations in the Speiser class, \rm
Preprint 2004.
\bibitem[CY]{cy} W. Cherry, Z. Ye, \em Nevanlinna's Theory of Value Distribution, \rm
Spinger Monographs in Mathematics (2001).
\bibitem[DU1]{du1} M. Denker, M. Urba\'nski, \em On the existence of conformal measures,  \rm
Trans. A.M.S. 328 (1991), 563-587.
\bibitem[DU2]{du2} M. Denker, M. Urba\'nski, \em Ergodic theory of
equilibrium states, \rm Nonlinearity 4 (1991), 103-134.
\bibitem[Er1]{eremenko VW} A. Eremenko \em On the iterations of entire functions, \rm
Dynamical systems and the ergodic theory, Banach Center Publication 23 (1989), 339-345.
\bibitem[Er2]{eremenko} A. Eremenko \em Ahlfors' contribution to the theory of
meromorphic functions  \rm .
\bibitem[EL]{el} A.E. Eremenko, M.Yu. Lyubich,  \em Dynamical properties of some classes of entire functions,   \rm
Ann. Inst. Fourier, Grenoble 42, 4 (1992), 989-1020.
\bibitem[H]{hille} E. Hille \em Analytic function theory, Vol. II, \rm Ginn (1962).
\bibitem[Hk]{hinkkanen} A. Hinkkanen, \em A sharp form of Nevanlinna's second fundamental theorem,
\rm Invent. Math. 108 (1992), 549-574.
\bibitem[IM]{im} C. Ionescu-Tulcea, G. Marinescu, \em Th\'eorie ergodique pour des classes d'operations
non-compl\`etement continues,   \rm Ann. Math. 52, (1950), 140-147.
\bibitem[Iv]{iv} F. Iversen \em Recherches sur les fonctions inverses dess fonctions m\'eromorphes,
\rm Th\`ese de Helsingfors (1914).
\bibitem[Ka]{ka} T. Kato, \em Perturbation theory for linear operators,   \rm Springer (1995).
\bibitem[KU1]{ku1} J. Kotus, M. Urba\'nski, \em Conformal, Geometric and invariant measures
for transcendental expanding functions, \rm Math. Annalen. 324 (2002),
619-656.
\bibitem[KU2]{ku2} J. Kotus, M. Urba\'nski, \em Geometry and ergodic theory of non-recurrent
elliptic functions, J. d'Analyse Math. 93 (2004), 35-102. \rm
\bibitem[KU3]{ku3} J. Kotus, M. Urba\'nski, \em The dynamics and geometry of the Fatou functions,
 \rm  Discrete \& Continuous Dyn. Sys. 13 (2005), 291-338.
\bibitem[KU4]{kuap} J. Kotus, M. Urba\'nski, \em Fractal Measures and Ergodic Theory
of Transcendental Meromorphic Functions, \rm Preprint 2004.
\bibitem[L1]{l1} M. Yu. Lyubich, \em Some typical properties of the dynamics of rational maps,   \rm
Russian Math. Surveys, 8, 5 (1983) 154-155.
\bibitem[L2]{l2} M. Yu. Lyubich, \em The dynamics of rational transforms: the topological picture,   \rm
Russian Math. Surveys, 41, 4 (1986) 43-117.
\bibitem[MSS]{mss} R. Ma\~{n}\'e, P. Sad, D. Sullivan, \em On the dynamics of rational maps,   \rm
Ann. Scient. Ec. Norm. Sup. 4e s\'erie, 16 (1983), 193-217.
\bibitem[McM]{mcm} C. McMullen, \em Area and Hausdorff dimension of Julia set of entire
functions, \rm Trans. A.M.S. 300 (1987), 329-342.
\bibitem[Mat]{mat} P. Mattila, \em Geometry of sets and measures in euclidean
spaces, \rm Cambridge Studies in Advanced Mathematics 44, Cambridge
University Press, 1995.
\bibitem[MdU]{mdu} D. Mauldin and M. Urba\'nski, \em Dimensions and measures
in infinite iterated function systems, \rm Proc. London Math. Soc.
(3) 73 (1996) 105-154.
\bibitem[My1]{my2} V. Mayer, \em Rational functions without conformal
measures on the conical set, \rm Preprint 2002.
\bibitem[My2]{my} V. Mayer, \em The size of the Julia set of
meromorphic functions, \rm Preprint 2005.
\bibitem[MyU]{myu1} V. Mayer, M. Urba\'nski, \em Gibbs and equilibrium
 measures for elliptic functions, \rm  Math. Zeitschrift 250 (2005), 657-683.
\bibitem[Nev1]{nev'} R. Nevanlinna, \em Eindeutige analytische
Funktionen, \rm Springer Verlag, Berlin (1953).
\bibitem[Nev2]{nev} R. Nevanlinna, \em Analytic functions, \rm Springer Verlag,
Berlin (1970).
\bibitem[Pa]{pat} S. J. Patterson, \em The limit set of a Fuchsian group, \rm
Acta Math. 136 (1976), 241-273.
\bibitem[PU]{pu} F. Przytycki, M. Urba\'nski, \em Fractals in the Plane - the
Ergodic Theory Methods, \rm available on Urba\'nski's webpage, to appear
Cambridge Univ. Press.
\bibitem[RS]{stri} P.J. Rippon, G.M. Stallard,  \em Iteration of a
class of hyperbolic meromorphic functions, \rm
Proc. of the AMS, Vol. 127, Nr. 11 (1999), 3251-3258.
\bibitem[R]{r} D. Ruelle, \em Repellers for real analytic maps,
\rm Ergod. Th. \& Dynam. Sys., 2 (1982), 99-107.

\bibitem[Sk]{sk} Z. Slodkowski, \em Holomorphic motions and polynomial hulls,   \rm
Proc. Amer. Math. Soc. 111 (1991), 347-355.

\bibitem[Su]{sul} D. Sullivan, \em Seminar on conformal and hyperbolic geometry,
\rm Preprint IHES  (1982).
\bibitem[UZ1]{uz0} M. Urba\'nski, A. Zdunik, \em The finer geometry and dynamics of
exponential family, \rm Michigan Math. J. 51 (2003), 227-250.
\bibitem[UZ2]{uz1}M. Urba\'nski, A. Zdunik, \em Real analyticity of
  Hausdorff dimension of
finer Julia sets of exponential family, \rm Ergod. Th. \& Dynam. Sys.
24 (2004), 279-315.
\bibitem[UZ3]{uzdnh} M. Urba\'nski, A. Zdunik, \em Geometry and
ergodic theory of non-hyperbolic exponential maps, \rm Preprint 2003,
to appear Trans. Amer. Math. Soc.
\end{thebibliography}
\end{document}